\documentclass[twocolumn]{IEEEtran}

\usepackage{amssymb}

\usepackage{setspace}

\usepackage{graphicx}
\usepackage{float}
\usepackage{graphicx}
\usepackage{caption}
\usepackage{amssymb}
\usepackage{amsmath}

\usepackage[boxed]{algorithm}
\usepackage{algorithmic}
\usepackage{longtable}
\usepackage{subfig}
\usepackage{epstopdf}
\usepackage{mathpazo}

\usepackage[boxed]{algorithm}
\usepackage{algorithmic}

\usepackage[hyphens]{url}

\usepackage{cite}

 \usepackage[boxed]{algorithm}
\usepackage{algorithmic}
\usepackage{color,xspace}
\usepackage{epstopdf}
\usepackage{balance}
\usepackage{bigints}
\usepackage{comment}
\usepackage[version=3]{mhchem}
\graphicspath{{Pictures/}}

\allowdisplaybreaks

\newtheorem{rmk}{Remark}

\newcommand{\soc}{{\rm SoC}\xspace}
\newcommand{\lqcwfts}{{\tt LQCwFTS}\xspace}
\newcommand{\lqt}{{\tt LQT}\xspace}
\newcommand{\sslqt}{{\tt SS-LQT}\xspace}
\newcommand{\bunderline}[1]{\underline{#1\mkern-4mu}\mkern4mu }

\newcommand{\T}{{\top}}

\begin{document}

\title{Health-Aware and User-Involved Battery Charging Management for Electric Vehicles: Linear Quadratic Strategies
\author{Huazhen Fang, Yebin Wang, and Jian Chen}
\thanks{H. Fang is with the Department of Mechanical Engineering,
University of Kansas,
Lawrence, KS 66045, USA (e-mail: fang@ku.edu).}
\thanks{Y. Wang is with the Mitsubishi Electric Research
Laboratories, Cambridge, MA 02139, USA (e-mail: yebinwang@ieee.org).}
\thanks{J. Chen is with the State Key Laboratory of Industrial Control Technology and College of Control
Science and Engineering, Zhejiang University, Hangzhou 310027, China
(e-mail: jchen@zju.edu.cn).}}

\markboth{Preprint Submitted to IEEE Transactions on Control Systems Technology}{}

\maketitle

\vspace{-15mm}

\begin{abstract} 
This work studies control-theory-enabled intelligent charging management for battery systems in electric vehicles (EVs). Charging is crucial for the battery performance and life as well as a contributory factor to a user's confidence in or anxiety about EVs. For existing practices and methods, many run with a lack of battery health awareness during charging, and none includes the user needs into the charging loop. To remedy such deficiencies, we propose to perform charging that, for the first time, allows the user to specify charging objectives and accomplish them through dynamic control, in addition to suppressing the charging-induced negative effects on battery health. Two charging strategies are developed using the linear quadratic control theory. Among them, one is based on control with fixed terminal charging state, and the other on tracking a reference charging path. They are computationally competitive, without requiring real-time constrained optimization needed in most charging techniques available in the literature. A simulation-based study demonstrates their effectiveness and potential. It is anticipated that charging with health awareness and user involvement guaranteed by the proposed strategies will bring major improvements to not only the battery longevity but also EV user satisfaction.
\end{abstract}


\IEEEpeerreviewmaketitle

\section{Introduction}

\IEEEPARstart{H}olding the promise for reduced fossil fuel use and air pollutant emissions, electrified transportation has been experiencing a surge of interest in recent years. Over 330,000 plug-in electric vehicles (EVs) are on the road in the United States as of May 2015~\cite{EDTA:2014}, with strong growth foreseeable in the coming decades. Most EVs rely on battery-based energy storage systems, which are crucial for the overall EV performance as well as consumer acceptance. Associated with this trend, the past years have witnessed a growing body of work on battery management research, e.g., state-of-charge (SoC) estimation to infer the amount of energy available in a battery, state-of-health
(SoH) estimation to track the battery's aging status, thermal monitoring to avoid abnormal heat buildup~\cite{Domenico:DSMC:2010,Fang:CEP:2014,Wang:TCST:2015,Kim:JPS:2015,Fang:JPS:2014,Smith:TCST:2010,Bartlett:TCST:2016,Dey:TCST:2015,Moura:DSMC:2013,Lin:TCST:2013}. Another essential yet less explored problem in the battery use is the charging strategies. Improper charging, e.g., charging with a high voltage or current density, can induce the rapid buildup of internal stress and resistance, crystallization and other negative effects~\cite{Suthar:PCCP:2013,Spotnitz:JPS:2003,Bergveld:Springer:2002,Catherino:JPS:2004}. The consequence is fast capacity fade and shortened life
cycle, and even safety hazards in the extreme case, eventually impairing the consumer confidence.

\begin{figure}  \centering
    \includegraphics[height=0.2\textwidth]{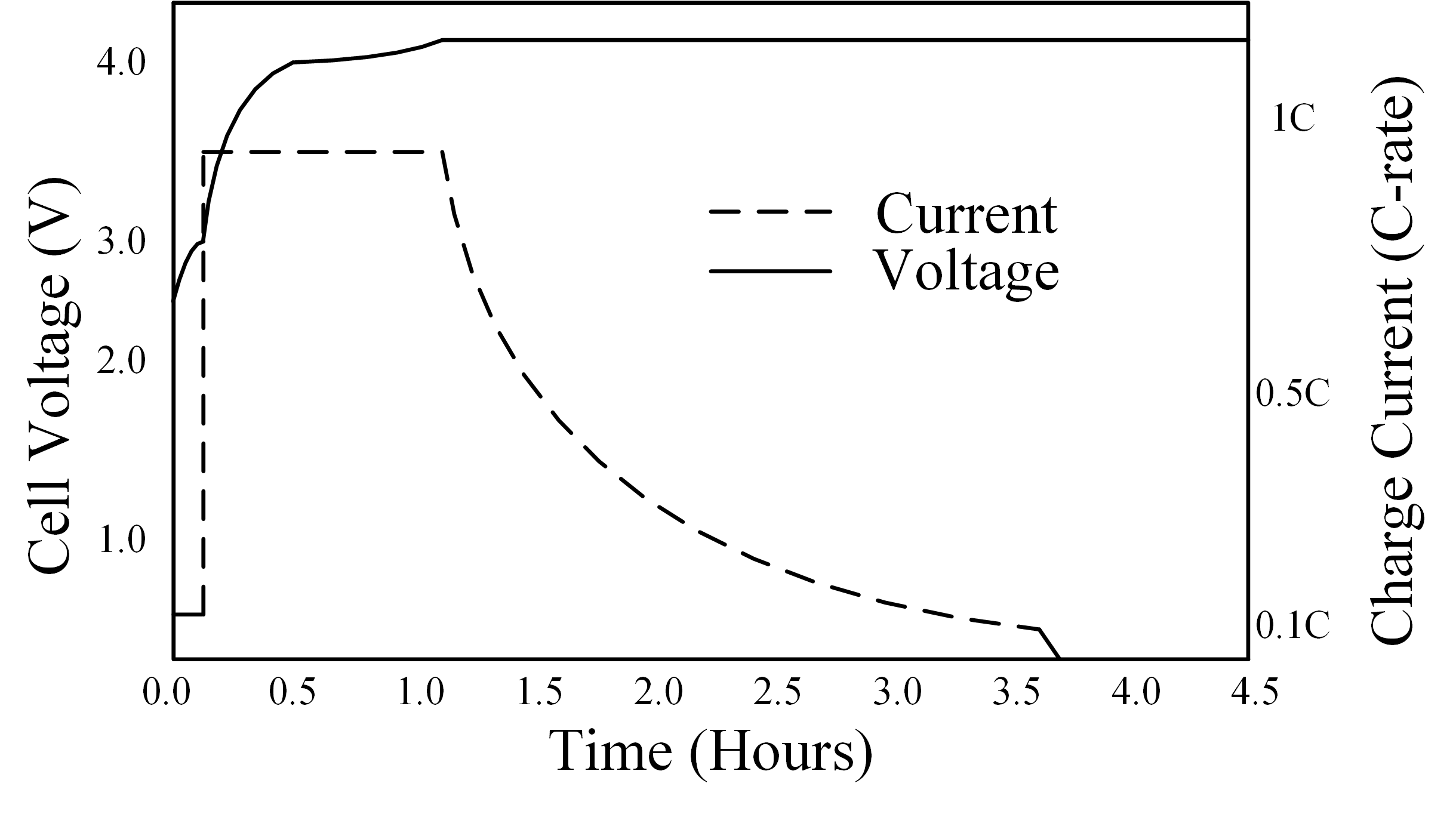}
  \caption{ Constant-current/constant-voltage charging.}
\label{CCCV-Charging}
\end{figure}

\begin{figure}  \centering
    \includegraphics[height=0.2\textwidth]{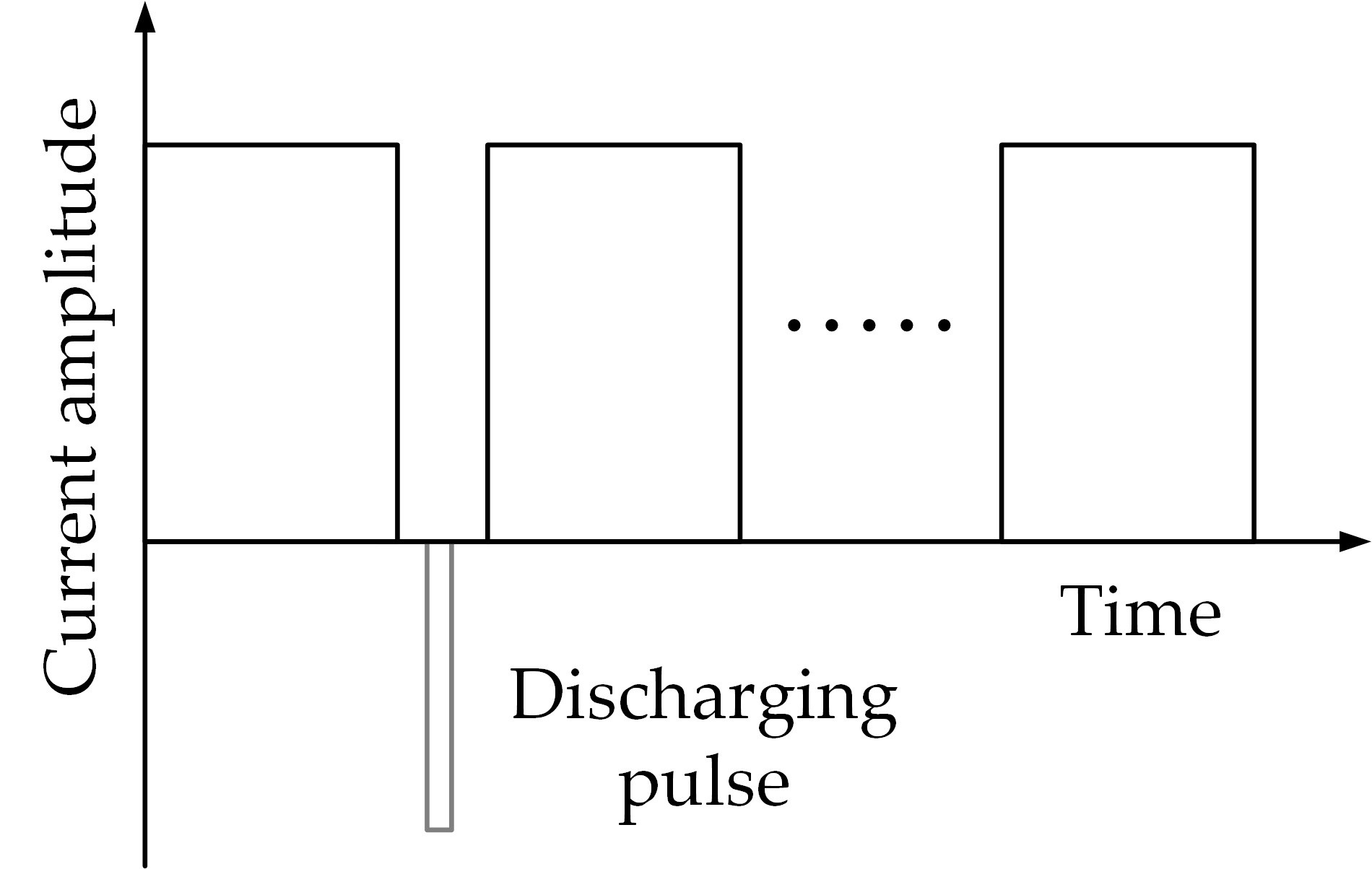}
  \caption{ Pulse charging and burp charging.}
\label{Pulse-Charging}
\end{figure}

{\em Literature review:}
The popular charging ways, especially for inexpensive lead-acid batteries used for cars and backup power systems, are to apply a constant voltage or force a constant current flow through the battery~\cite{Young:Springer:2012}. Such methods, though easy to implement, can lead to serious detrimental effects for the battery.
One improvement is the constant-current/constant-voltage charging~\cite{Rahn:Wiley:2013,Young:Springer:2012}, which is illustrated in Figure~\ref{CCCV-Charging}. Initially, a trickle charge (0.1C or even smaller) is used for depleted cells, which produces a rise of the voltage. Then a constant current between 0.2C and 1C is applied. This stage ends when the voltage increases to a desired level. The mode then switches to constant voltage, giving a diminishing current to charge. Yet the implementation is empirical here, with the optimal determination of the charge regimes remaining in question~\cite{Wong:JPS:2008}. In recent years, pulse charging has gained much interest among practitioners. Its current profile is based on pulses, as shown in Figure~\ref{Pulse-Charging}. Between two consecutive pulses is a short rest period, which allows the electrochemical reactions to stabilize by equalizing throughout the bulk of the electrode before the next charging begins. This brief relaxation can accelerate the charging process, reduce the gas reaction, inhibit dendrite growth and slow the capacity fade~\cite{Lam:JPS:1995,Purushothaman:JES:2006,Aryanfar:JPCL:2014}. Its modified
version, burp charging, applies a very short negative
pulse for discharging during the rest period , see Figure~\ref{Pulse-Charging}, in order
to remove the gas bubbles that have appeared the electrodes.

A main issue with the above methods is the lack of an effective feedback-based regulation mechanism. With an open-loop architecture, they simply take energy from power supply and put it into the battery. As a result, both the charging dynamics and the battery's internal state are not well exploited to control the charging process for better efficiency and health protection. This motivates the deployment of closed-loop and model-based control. Constrained optimal control is applied in~\cite{Klein:ACC:2011,Suthar:PCCP:2013,Yan:Energies:2011,Perez:Mechatronics:2015}, in conjunction with electrochemical or equivalent circuit models, to address fast charging subject to input, state and temperature constraints.  {In this direction, fast constrained optimization has been leveraged recently in~\cite{Torchio:ACC:2015,Liu:DSCC:2015} to reduce the computational cost and push forward real-time charging control.}
With the ability of dealing with uncertain parameters, adaptive control is used for energy-efficient fast charging in~\cite{Wai:IET-PE:2012}.
Based on the Pontryagin minimum principle, optimal control design of charging/discharging is studied in~\cite{Wang:AUTO:2012} to maximize the work that a battery can perform over a given duration while maintaining a desired final energy level. However, we observe that the research effort for feedback-controlled charging has remained limited to date. The existing works are mostly concerned with the fast charging scenario and employ a restricted number of investigation tools, thus presenting much scope for further work. 


{\em Research motivation:} In this paper, we propose to perform control-based EV charging management in a {\em health-aware} and {\em user-involved} way. Since the battery system is the heart as well as the most expensive component of an EV, health protection during charging is of remarkable importance to prevent performance and longevity degradation. As such, it has been a major design consideration in the controlled charging literature mentioned above. Furthermore, we put forward that the user involvement, entirely out of consideration in the state of the art, will bring significant improvements to charging. Two advantages at least will be created if the user can give the charging management system some commands or advisement about the charging objectives based on his/her immediate situation. The first one will be improved battery health protection against charging-induced harm. Consider two scenarios: 1) after arriving at the work place in the morning, a user leaves the car charging at the parking point with a forecast in mind that the next drive will be in four hours; 2) he/she will have a drive to the airport in one hour, and a half full capacity will be enough. In both scenarios, the user needs can be translated into charging objectives (e.g., charge duration and target capacity). The charger then can make wiser, more health-oriented charging decisions when aiming to meet the user specifications with such information, rather than pumping, effectively but detrimentally, the maximum amount of energy into the batteries within the minimum duration. Second, a direct and positive impact on user satisfaction will result arguably, because offering a user options to meet his/her varying and immediate charging needs not only indicates a better service quality, but also enhances his/her perception of level of involvement.

{\em Statement of contributions:} We will build health-aware and user-involved charging strategies via exploring two problems. The first one is {\em charging with fixed terminal charging state}. In this case, the user will give target SoC and charging duration, which will be incorporated as terminal state constraint. The second problem is {\em tracking-based charging}, where the charging is implemented via tracking a charge trajectory. The trajectory is generated on the basis of user-specified objectives and battery health conditions. The solutions, developed in the framework of {\em linear quadratic optimal control}, will be presented as controlled charging laws expressed in explicit equations. The proposed methods differ from those in the literature, e.g.,~\cite{Klein:ACC:2011,Suthar:PCCP:2013,Yan:Energies:2011,Wai:IET-PE:2012,Perez:Mechatronics:2015,Wang:AUTO:2012} in either of both of the following two aspects: 1) from the viewpoint of application, they keep into account both user specifications and battery health --- such a notion is unavailable before and will have a potential impact on improving the existing charging practices; 2) technically, they, though based on optimization of quadratic cost functions, do not require real-time constrained optimization needed in many existing techniques~\cite{Klein:ACC:2011,Suthar:PCCP:2013,Yan:Energies:2011,Perez:Mechatronics:2015} and thus are computationally more attractive. In addition, the linear quadratic control is a fruitful area, so future expansion of this work can be aided with many established results and new progresses, e.g.,~\cite{Lewis:Wiley:2012,Duncan:IEEE-TAC:1999,Duncan:IEEE-TAC:2013,Lee:AUTO:2000}.

{\em Organization:} The rest of the paper is organized as follows. Section~\ref{Model} introduces an equivalent circuit model oriented toward describing the battery charging dynamics. Section~\ref{Charging-startegies} presents the development of charging strategies. Section~\ref{Charging-with-fixed-terminal-state} studies the charging with fixed terminal charging state specified by the user. In Section~\ref{Tracking-based-charging}, tracking-based charging is investigated. Section~\ref{Simulation} offers numerical results to illustrate the effectiveness of the design. Finally, concluding remarks are gathered in Section~\ref{Conclusions}.

\begin{figure*}[htbp]
\normalsize
\begin{equation}\label{RC-model-Eq}
\left\{
\begin{aligned}
\left[\begin{matrix}\dot Q_b(t) \cr \dot Q_s(t) \end{matrix}\right]&=
\left[\begin{matrix} -\dfrac{1}{C_b(R_b+R_s)} & \dfrac{1}{C_s(R_b+R_s)} \cr \dfrac{1}{C_b(R_b+R_s)} & -\dfrac{1}{C_s(R_b+R_s)} \end{matrix}\right] \left[\begin{matrix} Q_b(t) \cr  Q_s(t) \end{matrix}\right] + \left[\begin{matrix} \dfrac{R_s}{R_b+R_s} \cr \dfrac{R_b}{R_b+R_s} \end{matrix}\right] I(t)\\
V(t)&= \left[ \begin{matrix} \dfrac{R_s}{C_b(R_b+R_s)} &  \dfrac{R_b}{C_s(R_b+R_s)}\end{matrix} \right] \left[\begin{matrix} Q_b(t) \cr  Q_s(t) \end{matrix}\right] + \left( R_o+ \frac{R_bR_s}{R_b+R_s}\right) I(t)
\end{aligned}.
\right.
\end{equation}
\hrulefill
\end{figure*}

\section{Charging Model Description}\label{Model}

\begin{figure}[t]  \centering
    \includegraphics[height=0.2\textwidth]{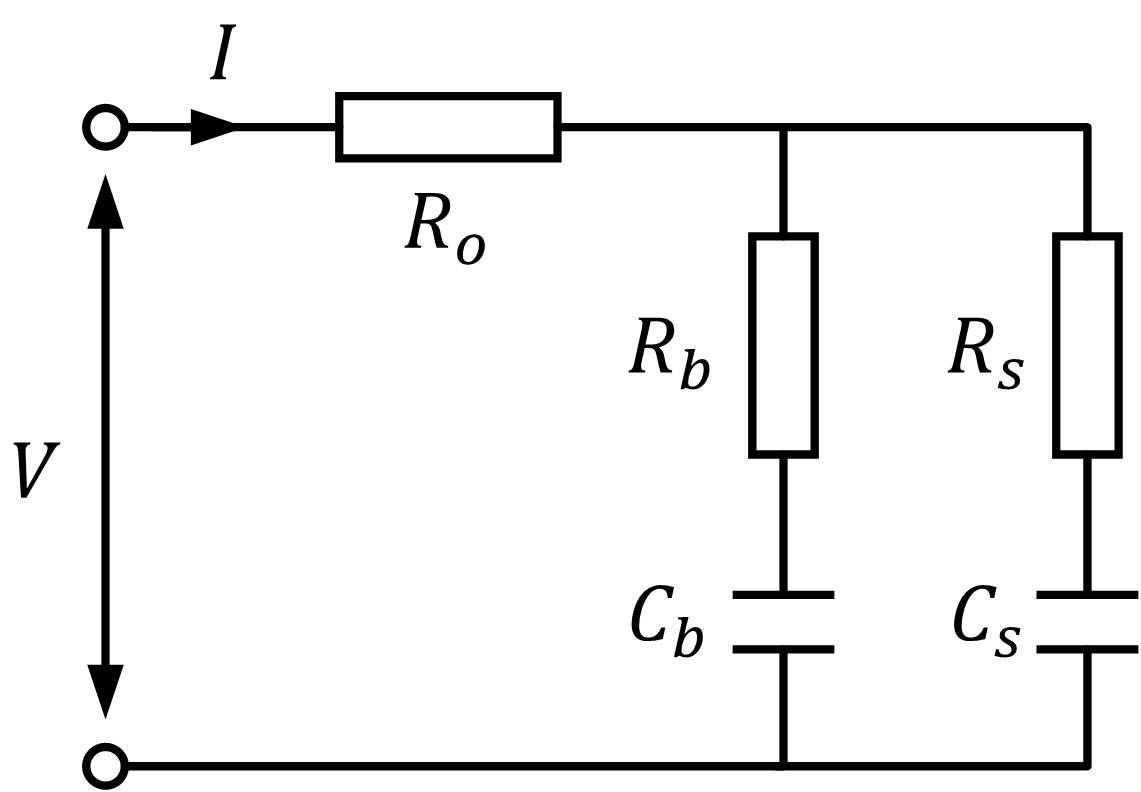}
  \caption{The battery RC model, where $R_o$, $R_s$-$C_s$ and $R_b$-$C_b$, respectively, simulates the resistance of the electrolyte, the surface region and the bulk inner part of an electrode.}
\label{RCM-1}
\end{figure}

While the energy storage within a battery results from complex electrochemical and physical processes, it has been useful to draw an analogy between the battery electrical properties and an equivalent circuit which consists of multiple linear passive elements such as resistors, capacitors, inductors and virtual voltage sources. While plenty of equivalent circuit models have been proposed, we focus our attention throughout the paper on a second-order resistance-capacitance (RC) model shown in Figure~\ref{RCM-1}.


Developed by Saft Batteries, Inc., this model was intended for the simulation of battery packs in hybrid EVs~\cite{Johnson:EVS:2000,Johnson:JPS:2002}. Identification of its parameters is discussed in~\cite{Sitterly:IEEE-TSE:2011}. 
As shown in Figure~\ref{RCM-1}, it consists of two capacitors and three resistors. {The resistor $R_o$ represents the electrolytic resistance within a battery cell. The double RC circuits in parallel are meant to simulate the migration of the electric charge during the charging (or discharging) processes. Specifically, the $R_s$-$C_s$ circuit accounts for the electrode surface region, which is exposed to the electrode-electrolyte interface; the $R_b$-$C_b$ circuit represents the bulk inner part of the electrode. 
Seeing a fast-speed transfer of the electric charge, the electrode surface is responsible for the high-frequency behavior during the charging processes and associated with the immediate amount of charge the battery can absorb. It, however, has a rather limited storage capacity. By contrast, the bulk electrode is where the majority of the electric charge is stored in chemical form. Since the diffusion of ions within the electrode proceeds at a relatively slower speed, the $R_b$-$C_b$ circuit makes up the low-frequency part of the charging response. This implies that $R_b \gg R_s$ and $C_b \gg C_s$. Note that it has been a significant effort to use the RC circuits to approximate electrochemical processes at different scales of time and frequency~\cite{Seaman:JPS:2014,Yuan:Springer:2010}. The state-space representation of the model is shown in~\eqref{RC-model-Eq}. It can be verified that this system is controllable and observable, indicating the feasibility of both controlled charging and state monitoring. }

Based on the model, the overall SoC is given by
\begin{align}\label{SoC-1}
\soc = \frac{Q_b-\bunderline{Q}_b + Q_s-\bunderline{Q}_s}{\bar Q_b-\bunderline{Q}_b+ \bar{Q}_s-\bunderline{Q}_s},
\end{align}
where $\bunderline{Q}_j$ and $\bar{Q}_j$ for $j=b,s$ denote the minimum and the maximum allowed charge held by the capacitor $C_j$, which represent the operating limits of the battery. When the equilibrium $V_b=V_s$ is reached, the SoC can be simply expressed as the linear combination of $\soc_b$ and $\soc_s$, i.e.,
\begin{align}\label{SoC-2}
\soc = \frac{C_b}{C_b+C_s} \soc_b +\frac{C_s}{C_b+C_s} \soc_s.
\end{align}


\begin{figure*}[t]  \centering
    \includegraphics[width=0.9\textwidth]{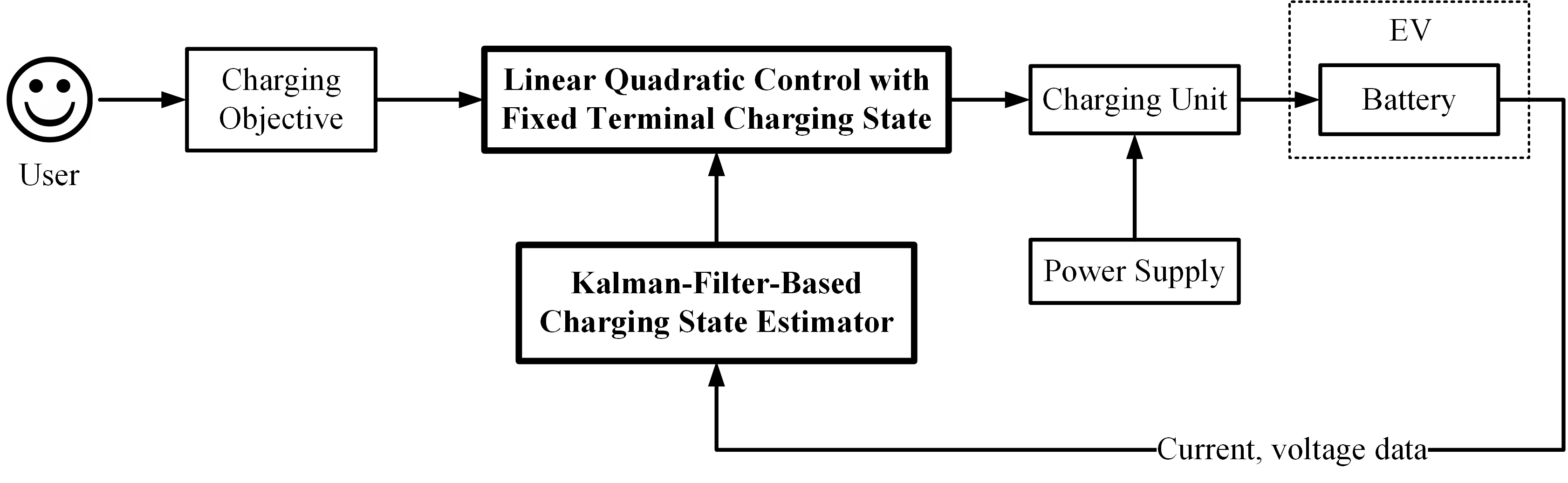}
  \caption{The schematic diagram for charging based on linear quadratic control with fixed terminal charging state.}
  \label{LQCwFTS-Diagram}
\end{figure*}

The RC model can well grasp the ``rate capacity effect'', which means that the total charge absorbed by a battery goes down with the increase in charging current as is often stated as the Peukert's law. To see this, consider that a positive current is applied for charging. Then both $Q_b$ and $Q_s$, and their voltages, $V_b$ and $V_s$, will grow. However, $V_s$ increases at a rate faster than $V_b$. When the current $I$ is large, the terminal voltage $V$, which is largely dependent on the fast increasing $V_s$, will grow quickly as a result. Then $V$ will reach the cut-off threshold in a short time. This will have the charging process terminated, though $Q_b$ still remains at a low level. Another essential phenomenon that can be well
approximated by this RC model is the ``recovery effect''.
upon an interruption of charging. That is, when the charging stops, the terminal voltage $V$ will see a transient decrease due to the charge transfer from $C_s$ to $C_b$.

To develop a digitally controlled charging scheme, the model in~\eqref{RC-model-Eq} is discretized with a sampling period of $t_s$. The discrete-time model takes the following standard form:
\begin{equation}\label{dsys}
\left\{
\begin{aligned}
x_{k+1} &=  A x_k + B u_k\\
y_k &= C x_k + D u_k
\end{aligned}.
\right.
\end{equation}
where $x=\left[ \begin{matrix} Q_b & Q_s \end{matrix} \right]^\T$, $u=I$, $y=V$, and $A$, $B$, $C$ and $D$ can be decided according to the discretization method applied to~\eqref{RC-model-Eq}.

{ Despite being linear and straightforward, the above RC model can satisfy the practical needs in many applications. This is because battery systems, e.g., those in electric vehicles, need to limit the minimum and maximum \soc during operation~\cite{Knutsen:2013,Markel:AABC:2006} for the purposes of safety, life, and a consistent power capability. Within this favorable \soc range, the battery behavior can be approximated as linear due to battery open-circuit-voltage profiles.}

For health consideration, we need to constrain the difference between voltages across $C_s$ and $C_b$, denoted as $V_s$ and $V_b$ respectively, throughout the charging process. Note that $V_s -V_b$ is the force that drives the migration of the charge from $C_s$ to $C_b$ during charging. It shares great resemblance with the gradient of the concentration of Li ions within the electrode created during charging causing the diffusion of ions. { This observation is further investigated in the Appendix. where the analogy between the voltage difference and the Li-ion concentration gradient is validated through a proof approximate equivalence between the model in~\eqref{RC-model-Eq} and the well-known single particle model (SPM) under certain conditions.} Too large a gradient value will cause internal stress increase, heating, solid-electrolyte interphase (SEI) formation and other negative side effects~\cite{Pinsona:JES:2013,Woodford:MIT:2013,Bandhauera:JES:2011}.
Mechanical degradation in the electrode and capacity fade will consequently happen. Thus to reduce the battery health risk, non-uniformity of the ion concentration should be suppressed during charging, {and this implies the necessity of suppressing the voltage difference}. It is also noteworthy that such a restriction should be implemented more strictly as the SoC increases, because the adverse effects of a large concentration difference  on the battery would be stronger then. 

Next, we will build the charging strategies on the basis of the above RC model. The development will be laid out in the framework of linear quadratic control, taking into account both health awareness and user needs.

\section{Health-Aware and User-Involved Charging Strategies} \label{Charging-startegies}

In this section, we will develop two charging strategies. For both, the user specifies the desired charging duration and target capacity.  The first strategy accomplishes the task via a treatment based on linear quadratic control subject to fixed terminal state resulting from the user objective. In the second case, charging is managed via tracking a charging trajectory which is produced from the user objective. A discussion of the strategies will follow.

\subsection{Charging with Fixed Terminal Charging State}\label{Charging-with-fixed-terminal-state}

A charging scenario that frequently arises is:  according to the next drive need, a user will inform the charging management system of his/her objective in terms of target SoC and charging duration. This can occur for overnight parking at home and daytime parking at the workplace, or when a drive to some place will set off in a predictable time. As afore discussed, the objective offered by the user, if incorporated into the dynamic charging decision making process, would create benefits for health  protection compared to fast charging.
This motivates us to propose a control-enabled charging system illustrated in Figure~\ref{LQCwFTS-Diagram}.
The charging objective given by the user is taken and and translated into the desired terminal charging state. A linear quadratic controller will compute online the charging current to apply so as to achieve the target state when the charging ends. Meanwhile, a charging state estimator will  monitor the battery status using the current and voltage measurements, and feed the information to the controller. In the following, we will present how to realize the above charging control.

From the perspective of control design, the considered charging task can be formulated as an optimal control problem, which minimizes a cost function commensurate with the harm to health and subject to the user's goal. With the model in~\eqref{dsys}, the following linear quadratic control problem will be of interest:
\begin{equation}
\begin{aligned}\label{Charging-with-FFS}
\min_{u_0, u_1,\cdots, u_{N-1}} & \ \frac{1}{2} x_N^\T S_N x_N \\ & \hspace{-5mm} + \frac{1}{2} \sum_{k=0}^{N-1} \left(x_k^\T G^\T Q_k G x_k  + u_k^\T R u_k \right), \\
\mbox{subject to} & \ x_{k+1} = A x_k + Bu_k, \ x_0\\
& \ x_N = \bar x.
\end{aligned}
\end{equation}
where $S_N\geq 0$, $Q_k\geq 0$, $R>0$ and 
\begin{align*}
G = \left[ \begin{matrix}
-\frac{1}{C_b} &  \frac{1}{C_s}
\end{matrix}
\right].
\end{align*}
In above, $G x_k$ represents the potential difference between $C_b$ and $C_s$, and the time range $N$ and the final state $\bar x$ are generated from the user-specified charging duration and target \soc. Note that the battery should be at the equilibrium point with $V_b = V_s$ in the final state and that using~\eqref{SoC-1}-\eqref{SoC-2}, $\bar x$ can be easily determined from the specified SoC value. The quadratic cost function thus intends to constrain the potential difference and magnitude of the charging current during the charging process. The minimization is subject to both the state equation and the fixed terminal state. The weight coefficient $Q_k$ should be chosen in a way such that it increases over time, in order to offer stronger health protection that is needed as the SoC builds up. {It is noted that $x_N^\T S_N x_N$ represents a general formulation of the terminal cost, to which different options can be assigned. It vanishes, for example, if $S_N = 0$. Or it can be set as $S_N = G^\T Q_N G$ to constrain the voltage difference in the end state. However, since subjected to the hard constraint $x_N = \bar x$, the end state would reach the desired point regardless of $S_N$.}

Resolving the problem in~\eqref{Charging-with-FFS} will lead to a state-feedback-based charging strategy, which can be expressed in a closed-form~\cite{Lewis:Wiley:2012}:
\begin{align}\label{LQwFTS-1}
K_k &= (B^\T S_N B +R)^{-1}B^\T S_{k+1}A,\\ \label{LQwFTS-2}
S_k &= A^\T S_{k+1} (A-BK_k)+Q_k,\\ \label{LQwFTS-3}
T_k &= (A-BK_k)^\T T_{k+1}, \ T_N = I,\\ \label{LQwFTS-4} \nonumber
P_k &= P_{k+1} - T_{k+1}^\T B (B^\T S_{k+1} B +R)^{-1} B^\T T_{k+1}, \\  P_N & = 0, \\ \label{LQwFTS-5}
K^u_k &= \left(B^\T S_{k+1} B +R\right)^{-1} B^\T,\\ \label{SF-FTS-Control-Law} \nonumber
u_k &= - \left(K_k - K_k^u T_{k+1} P_k^{-1} T_k^\T \right) x_k \\& \quad \quad- K_k^u T_{k+1} P_k^{-1} \bar x.
\end{align}
This procedure comprises offline backward computation of the matrices $K_k$, $S_k$, $T_k$, $P_k$ and $K_k^u$ from the terminal state and online forward computation of the control input (i.e., charging current) $u_k$.

The state variable $x_k$ is not measurable directly in practice, so its real-world application necessitates the conversion of the above state-feedback-based strategy to an output-feedback-based one. One straightforward avenue to achieve this is to replace $x_k$ by its prediction $\hat x_k$. This is justifiable by the certainty equivalence principle, which allows the optimal output-feedback control design to be divided into the separate designs of an optimal state-feedback control and an optimal estimator~\cite{Bryson:Taylor:1975}. The optimal estimation can be treated via minimizing 
\begin{align}\label{KF-Cost-Func} \nonumber
\min_{x_0,x_1,\cdots, x_k} & \  \frac{1}{2} (x_0 - \hat x_0)^\T \Sigma_0^{-1} (x_0 - \hat x_0) \\ &  + \frac{1}{2} \sum_{i=0}^{k-1} w_i^\T {\Pi}^{-1} w_i + \frac{1}{2} \sum_{i=0}^{k} v_i^\T {\Lambda}^{-1} v_i ,
\end{align}
where $\Sigma_0>0$, $\Pi>0$, $\Lambda>0$, and
\begin{align*}
w_k &= x_{k+1}-Ax_k-Bu_k,\\
v_k & = y_k - Cx_k-D u_k.
\end{align*}
The one-step-forward Kalman predictor will result from solving~\eqref{KF-Cost-Func}, which is given by
\begin{align} \label{KF-1}
L_k &=  A\Sigma_k C^\T (C\Sigma_k C^\T +\Lambda)^{-1},\\ \label{KF-2}
\hat x_{k+1} &= A \hat x_k + B u_k +L_k (y_k - C \hat x_k -D u_k),\\ \label{KF-3} \nonumber
\Sigma_{k+1} &= A \Sigma_k A^\T +\Pi - A\Sigma_k C^\T \\ & \quad\quad \cdot (C\Sigma_k C^\T +\Lambda)^{-1}C\Sigma_kA^\T.
\end{align}
Substituting $x_k$ with its estimate $\hat x_k$, the optimal control law in~\eqref{SF-FTS-Control-Law} will become
\begin{align} \label{LQwFTS-6}
u_k &= - \left(K_k - K_k^u T_{k+1} P_k^{-1} T_k^\T \right) \hat x_k - K_k^u T_{k+1} P_k^{-1} \bar x.
\end{align}

\begin{table*}[!hbt]
\centering
\begin{tabular}{|l|}
\hline \\
Offline backward computation (from time $N$ to $0$)\\
\quad $K_k = (B^\T S_N B +R)^{-1}B^\T S_{k+1}A$\\ 
\quad $S_k = A^\T S_{k+1} (A-BK_k)+Q_k$\\
\quad $T_k = (A-BK_k)^\T T_{k+1}, \ T_N = I$\\
\quad $P_k = P_{k+1} - T_{k+1}^\T B (B^\T S_{k+1} B +R)^{-1} B^\T T_{k+1}, \ P_N = 0$\\
\quad $K^u_k = \left(B^\T S_{k+1} B +R\right)^{-1} B^\T$\\
\\
\hline
\\
Online forward computation (from time 0 to $N$)\\
\quad {\em Battery state prediction}\\
\quad \quad $L_k =  A\Sigma_k C^\T (C\Sigma_k C^\T +\Lambda)^{-1}$\\
\quad \quad$\hat x_{k+1} = A \hat x_k + B u_k +L_k (y_k - C \hat x_k -D u_k)$\\
\quad \quad$\Sigma_{k+1} = A \Sigma_k A^\T +\Pi - A\Sigma_k C^\T (C\Sigma_k C^\T +\Lambda)^{-1}C\Sigma_kA^\T$\\
\quad {\em Charging decision} \\
\quad \quad$u_k = - \left(K_k - K_k^u T_{k+1} P_k^{-1} T_k^\T \right) \hat x_k - K_k^u T_{k+1} P_k^{-1} \bar x$\\
\\
\hline
  \end{tabular}
  \caption{The \lqcwfts charging strategy (Linear Quadratic Control with Fixed Terminal State).}
  \label{LQCwFTS}
\end{table*}

\begin{figure*}[!t]  \centering
    \includegraphics[width=0.95\textwidth]{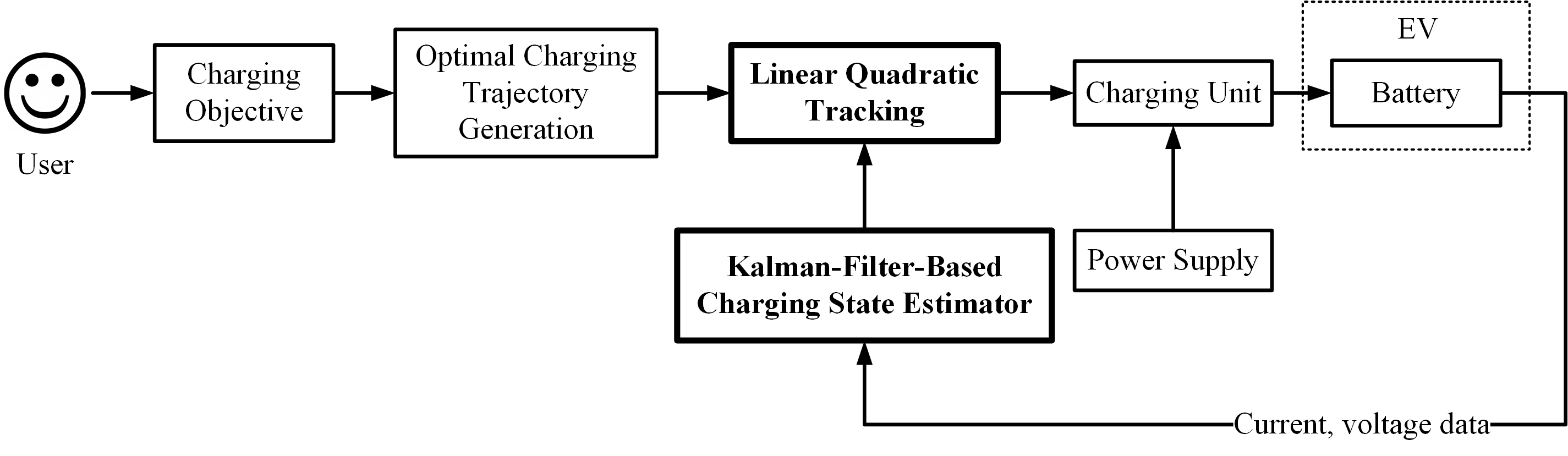}
  \caption{The schematic diagram for charging based on linear quadratic tracking.}
\label{LQT-Diagram}
\end{figure*}

Putting together~\eqref{LQwFTS-1}-\eqref{LQwFTS-5}, \eqref{KF-1}-\eqref{KF-3} and~\eqref{LQwFTS-6}, we will obtain a complete description of the charging method via linear quadratic control with fixed terminal state, which is named \lqcwfts and summarized in Table~\ref{LQCwFTS}. The \lqcwfts method performs state prediction at each time instant, and then feeds the predicted value, which is a timely update about the battery's internal state, to generate the control input to charge the battery. Much of the computation for \lqcwfts can be performed prior to the implementation of the control law. The sequences, $K_k$, $S_k$, $T_k$, $P_k$ and $K_k^u$ can be computed offline, and then $K_k$, $K_k^u T_{k+1} P_k^{-1} T_k^\T$ and $K_k^u T_{k+1} P_k^{-1}$ are stored for use when the control is applied. On the side of the Kalman prediction, offline computation and storage of $L_k$ can be done. Then the only work to do during charging is to compute the optimal state prediction and control input by~\eqref{KF-2} and~\eqref{LQwFTS-6}, reducing the computational burden.

\subsection{Charging Based on Tracking} \label{Tracking-based-charging}

Tracking-control-based charging is another way to guarantee health awareness and user objective satisfaction. A schematic of  its realization is shown in Figure~\ref{LQT-Diagram}. When a user specifies the charging objective, a charging trajectory can be generated. A charging controller will be in place to track the path. The trajectory generation will be conducted with a mix of prior knowledge of the battery electrochemistries, health awareness and user needs. It is arguably realistic that an EV manufacturer can embed trajectory generation algorithms into BMSs mounted on EVs, from which the user can select the one that best fits the needs when he/she intends to charge the EV. Leaving optimal charging trajectory generation for our future quest, we narrow our attention to the focus of path-tracking-based charging control here.

\begin{table*}
\centering
\begin{tabular}{|l|}
\hline \\
Offline backward computation (from time $N$ to $0$)\\
\quad $K_k  = (B^\T S_{k+1} B+R)^{-1} B^\T S_{k+1}A$\\
\quad $K_k^s = (B^\T S_{k+1} B+R)^{-1} B^\T$\\
\quad $S_k = A^\T S_{k+1} (A-BK_k) +Q $\\
\quad $s_k = (A-BK_k)^\T s_{k+1} + Q r_k,  \ s_N = S_N r_N $\\
\\
\hline
\\
Online forward computation (from time 0 to $N$)\\
\quad {\em Battery state prediction}\\
\quad \quad $L_k =  A\Sigma_k C^\T (C\Sigma_k C^\T +\Lambda )^{-1}$\\
\quad \quad$\hat x_{k+1} = A \hat x_k + B u_k +L_k (y_k - C \hat x_k -D u_k)$\\
\quad \quad$\Sigma_{k+1} = A \Sigma_k A^\T +\Pi - A\Sigma_k C^\T (C\Sigma_k C^\T +\Lambda)^{-1}C\Sigma_kA^\T$\\
\quad {\em Charging decision} \\
\quad \quad $u_k = -K_k \hat x_k +K_k^s s_{k+1}$\\
\\
\hline
  \end{tabular}
  \caption{The \lqt charging strategy (Linear Quadratic Tracking).}
  \label{LQT}
\end{table*}

Suppose that the user describes the target SoC and duration for charging, which are translated into the final state $\bar x$. Then a reference trajectory $r_k$ for $k=0,1,\cdots,N$ is calculated with $r_N = \bar x$. Note that the trajectory should constrain the difference between $V_b$ and $V_s$ to guarantee health. A linear quadratic state-feedback tracking can be considered for charging:
\begin{equation}
\begin{aligned} 
\min_{u_0,u_1,\cdots,u_{N-1}} &\ \frac{1}{2} \left(x_N-r_N\right)^\T S_N \left(x_N-r_N\right) \\ & \hspace{-16mm}  + \frac{1}{2} \sum_{k=0}^{N-1} \left[\left(x_k-r_k\right)^\T Q \left(x_k-r_k\right)  + u_k^\T R u_k \right],\\
{\rm subject \ to} & \quad  x_{k+1} = Ax_k+Bu_k, \ x_0
\end{aligned}
\end{equation}
where $S_N\geq 0$, $Q\geq 0$ and $R>0$.
Referring to~\cite{Lewis:Wiley:2012}, the optimal solution to the above problem is expressed as follows:
\begin{align}\label{LQT-1}
K_k  &= (B^\T S_{k+1} B+R)^{-1} B^\T S_{k+1}A,\\ \label{LQT-2}
K_k^s &= (B^\T S_{k+1} B+R)^{-1} B^\T,\\ \label{LQT-3}
S_k &= A^\T S_{k+1} (A-BK_k) +Q ,\\ \label{LQT-4}
s_k&= (A-BK_k)^\T s_{k+1} + Q r_k,  \ s_N = S_N r_N, \\ \label{LQT-Control-Law}
u_k &= -K_k x_k +K_k^s s_{k+1}.
\end{align}
Resembling~\eqref{LQwFTS-1}-\eqref{SF-FTS-Control-Law}, the execution of the above procedure is in a backward-forward manner. Specifically,~\eqref{LQT-1}-\eqref{LQT-4} are computed offline and backward prior to charging, and~\eqref{LQT-Control-Law} online and forward from the moment when charging begins.

Following lines analogous to the development of \lqcwfts, the output-feedback tracker for charging can be created based on~\eqref{LQT-1}-\eqref{LQT-Control-Law} running with the Kalman predictor in~\eqref{KF-1}-\eqref{KF-3}. That is,~\eqref{LQT-Control-Law} will use $\hat x_k$ rather than $x_k$ in practical implementation, i.e.,
\begin{align} \label{LQT-5}
u_k = -K_k \hat x_k +K_k^s s_{k+1}.
\end{align}

\begin{table*}
\centering
\begin{tabular}{|l|}
\hline \\
Offline computation of DAREs and gains\\
\quad $S = A^\T S A - A^\T S B (B^\T S B +R)^{-1} B^\T S A +Q$\\
\quad $\Sigma = A \Sigma A^\T - A \Sigma C^\T (C \Sigma C^\T + \Lambda )^{-1} C \Sigma A^\T +\Pi$\\
\quad $\bar K  = (B^\T S B+R)^{-1} B^\T S A$\\
\quad $\bar K^s = (B^\T S B+R)^{-1} B^\T$\\
\quad $\bar L =  A\Sigma C^\T (C\Sigma C^\T + \Lambda )^{-1}$
\\ 
\\
\hline
\\
Offline computation of $s_0$ (from time $N$ to 0)\\
\quad ${ s_k} = (A-B \bar K)^\T s_{k+1} + Q r_k,  \ s_N = S_N r_N $\\
\\
\hline
\\
Online forward computation (from time 0 to $N$)\\
\quad {\em Battery state prediction}\\
\quad \quad$\hat x_{k+1} = A \hat x_k + B u_k +\bar L (y_k - C \hat x_k -D u_k)$\\
\quad {\em Charging decision} \\
\quad \quad $s_{k+1} = (A-B\bar K)^{-\T} s_k - (A-B\bar K)^{-\T}  Q r_k$\\
\quad \quad $u_k = -\bar K \hat x_k + \bar K^s s_{k+1}$\\
\\
\hline
  \end{tabular}
  \caption{The \sslqt charging strategy (Steady-State Linear Quadratic Tracking).}
  \label{SS-LQT}
\end{table*}

Summarizing~\eqref{LQT-1}-\eqref{LQT-4},~\eqref{KF-1}-\eqref{KF-3} and~\eqref{LQT-5} will yield the linear quadratic tracking strategy, or \lqt, for charging, see Table~\ref{LQT}. Similar to the aforeproposed \lqcwfts, the \lqt can have much computation completed offline. Then only the Kalman state prediction and optimal tracking control~\eqref{LQT-5} need to be computed during the actual control run.

The computational cost of \lqt can be further reduced if we use its steady-state counterpart, making it more desirable in the charging application. The steady-state tracker is deduced as follows. It is known that, if $(A,B)$ is stabilizable and $(A,Q^{1\over 2})$ is detectable, $S_k$, as $N-k\rightarrow \infty$, will approach a unique stabilizing solution of the discrete algebraic Riccati equation (DARE)
\begin{align*}
S = A^\T S A - A^\T S B (B^\T S B +R)^{-1} B^\T S A +Q.
\end{align*}
Then $K_k$ and $K_k^s$ will approach their respective steady-state values, $\bar K$ and $\bar K^s$. In a similar way, the Kalman gain $L_k$ will achieve steady state $\bar L$ as $k\rightarrow \infty$ given the detectability of $(A,C)$ and stabilizability of $(A,Q^{1\over 2})$, which is the unique stabilizing solution to the DARE 
\begin{align*}
\Sigma = A \Sigma A^\T - A \Sigma C^\T (C \Sigma C^\T + \Lambda)^{-1} C \Sigma A^\T +\Pi.
\end{align*}
According to the DARE theory, $S$ and $\Sigma$ can be solved for analytically. With the steady-state gains $\bar K$, $\bar K^s$ and $\bar L$, the optimal prediction and control for charging will be
\begin{align}\label{SS-LQT-1}
u_k &= - \bar K \hat x_k +\bar K^s s_{k+1},\\ \label{SS-LQT-2}
\hat x_{k+1} &= A \hat x_k + B u_k +\bar L (y_k - C \hat x_k -D u_k).
\end{align}
If $(A-B\bar K)$ is invertible, the backward computation of $s_k$ can be substituted by the forward computation governed by
\begin{align} \label{SS-LQT-3}
s_{k+1} = (A-B\bar K)^{-\T} s_k - (A-B\bar K)^{-\T}  Q r_k.
\end{align}
Its implementation is initialized by $s_0$ computed offline by~\eqref{LQT-4}.
We refer to this suboptimal charging strategy~\eqref{SS-LQT-1}-\eqref{SS-LQT-3} as the steady-state \lqt, or \sslqt and outline it in Table~\ref{SS-LQT}. The \sslqt strategy, due to its exceptional simplicity, has more computational appeal in terms of time and space complexity.

\subsection{Discussion}

The following remarks summarize our discussion of the proposed charging strategies.

\begin{rmk}{\em (Soft-constraint-based health awareness).} As is seen, the proposed \lqcwfts, \lqt and \sslqt strategies incorporate the health awareness as part of the cost functions rather than hard constraints. This soft-constraint-based treatment will bring the primary benefit of computational efficiency and convenience. This compares with the techniques based on real-time constrained optimization, which are relatively more time-consuming and on occasions, face the issue that no feasible solution exists in the constrained region. {In the meantime, soft constraints are acknowledged as less powerful than hard constraints, e.g., ~\cite{Klein:ACC:2011,Suthar:PCCP:2013,Yan:Energies:2011,Perez:Mechatronics:2015,Torchio:ACC:2015,Liu:DSCC:2015}, in terms of preventing violation of certain physical limits during charging. However, we argue that the use of soft constraints
does not compromise the effectiveness of the proposed strategies to protect the battery health. This is fundmentally because the usual cause of an actual limit violation is too aggressive a charging current and an essential part of the proposed strategies is to suppress the magnitude of the charging current.  For instance, it is noted that the harm to health is associated with a weighted penalty for the \lqcwfts. When a proper weight $Q_k$ is selected, minimizing the penalty cost will ensure a sufficient consciousness of the health. }
\end{rmk}

\begin{rmk}{\em (Robustness of \sslqt).} The \sslqt strategy is based on a combination of linear quadratic tracker and a Kalman filter. Such a design may engender weak robustness in terms of gain and phase margins. To overcome this limitation, the loop transfer recovery can be used to build robust control design on the linear quadratic control structure~\cite{Lewis:Wiley:2012}.
\end{rmk}

{
\begin{rmk}{\em (Choice of $Q_k$ and $R$ for \lqcwfts).} When the weight coefficients $Q_k$ and $R$ take different values, the charging profiles generated by the \lqcwfts strategy will change accordingly. This implies the importance of finding appropriate $Q_k$ of $R$ for the implementation. A basic guideline is as follows:
\begin{itemize}

\item $Q_k$ should increase over time to suppress the use of a large current when \soc becomes larger, because of a battery's susceptibility increasing with SoC to the charging current. 

\item $Q_k\gg R$, because the $Q$-weighted term in $J$ is much smaller than the $R$-weighted term. 

\item The larger $Q_k$ is, the less aggressive charging action will be. However, the overall charging action also depends on the final state constraint.
\end{itemize}
\noindent It should be noted that the selection of $\alpha$ and $\beta$ is a multi-faceted issue, because it needs to account for both battery health protection and charging speed and more broadly, the economic cost and user satisfaction. Since these factors depend on specific application scenarios, we leave this issue for practitioners to resolve based on the above guideline.
\end{rmk}
}

\begin{rmk}{\em (Generality to other models).}
The proposed development has a potential applicability to other battery models. {First, the investigation, though based on a linear model, can be extended to nonlinear battery models. It is observed that, for variaous control-oriented battery models, the nonlinearity exists only in measurement equation that relates the state and applied current with the measured output voltage. Thus, an extension can be readily made by deploying a nonlinear Kalman filter for state estimation without changing the control structure.} We can also generalize the design to the well-known single particle model (SPM). This model represents each battery electrode as a spherical particle and delineates the migration of ions in and between the particles as a diffusion process~\cite{Chaturvedi:CSM:2010}. The PDE-based SPM can be converted into the standard linear state-space form, as shown in~\cite{Domenico:DSMC:2010}. Then following similar lines to this work, linear quadratic problems can be established and solved for charging tasks, where the difference of ion concentration gradients are constrained to penalize charging-induced harm. {It is also worth pointing out that extensions can be
made to accommodate the temperature dynamics as a means to
suppress the charging-induced heat build-up. Specifically for
the considered model in Equation (1), a thermal coupling can
be performed as shown in~\cite{Sitterly:IEEE-TSE:2011}. We can then follow similar
lines to accomplish the linear quadratic charging design based
on the modified model.}
\end{rmk}


\section{Numerical Illustration}\label{Simulation}

In this section, we present two simulation examples to illustrate the performance of the proposed charging strategies. 
Let us consider a lithium-ion battery described by the RC model in~\eqref{RC-model-Eq} with known parameters provided by Saft Inc. for hybrid electric vehicles, with $C_b = 82$ kF, $R_b=1.1 \ {\rm m}\Omega$, $C_s = 4.074$ kF, $R_s = 0.4 \ {\rm m}\Omega$, and $R_o = 1.2 \ {\rm m}\Omega$~\cite{Johnson:EVS:2000}. It has a nominal capacity of 7 Ah. The model is discretized by a sampling period of $t_s = 1$ s. The initial SoC is assumed to be $30\%$. The user will specify that certain SoC must be achieved within certain duration. 

\begin{figure*}\centering
    \subfloat[]{
    \includegraphics[width=0.45\textwidth]{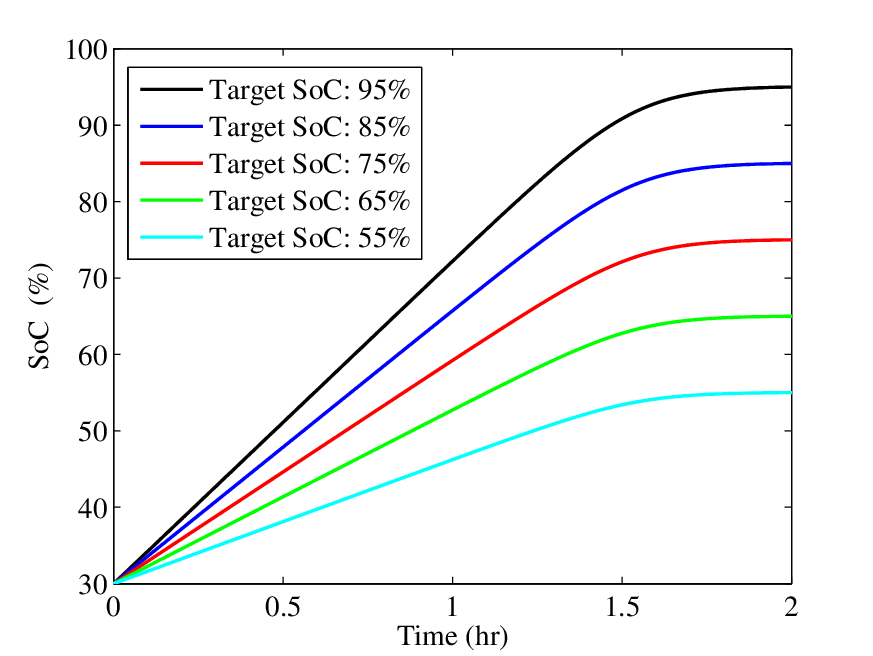}\label{Ex1-SoC}}
  \subfloat[]{
    \includegraphics[width=0.45\textwidth]{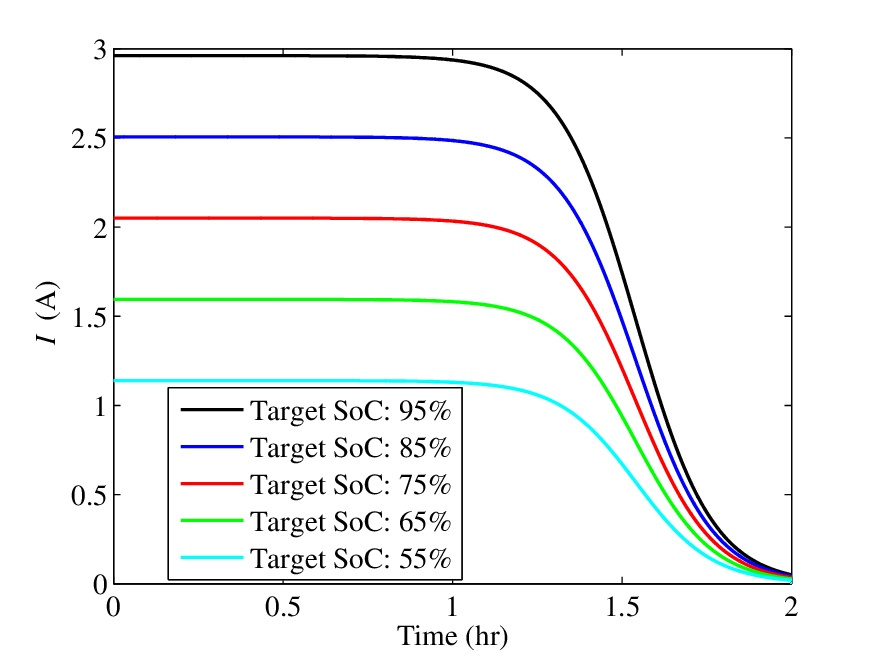}\label{Ex1-Current}}\\
	\subfloat[]{
   \includegraphics[width=0.45\textwidth]{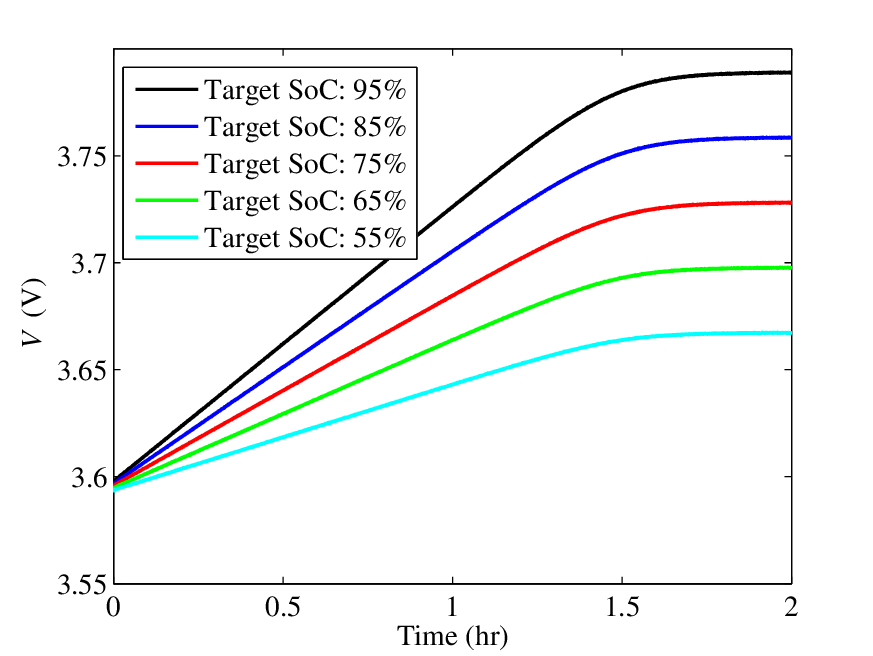}\label{Ex1-Volt-Output}}
   \subfloat[]{
   \includegraphics[width=0.45\textwidth]{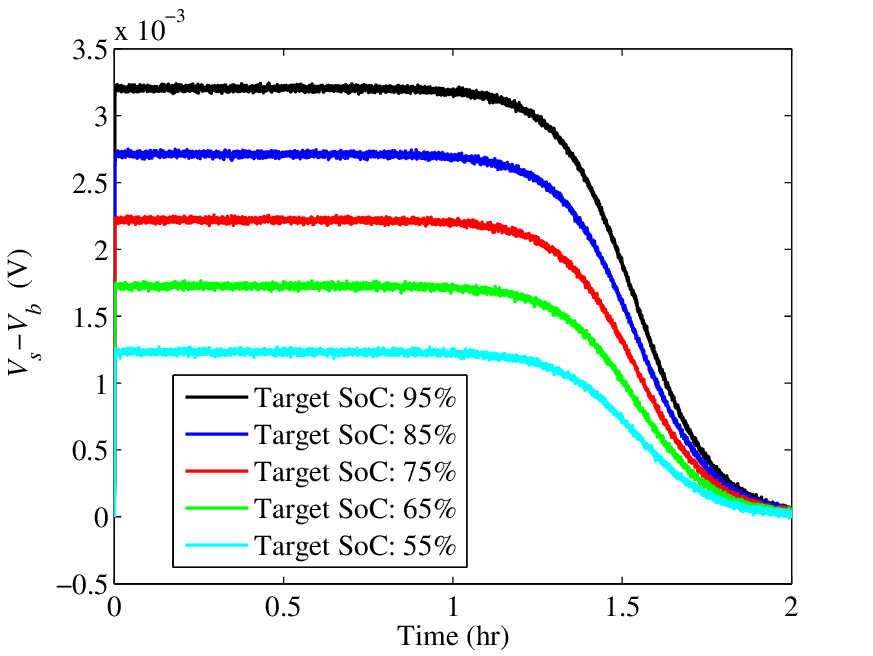}\label{Ex1-Voltage-Diff}}\\
\subfloat[]{
   \includegraphics[width=0.45\textwidth]{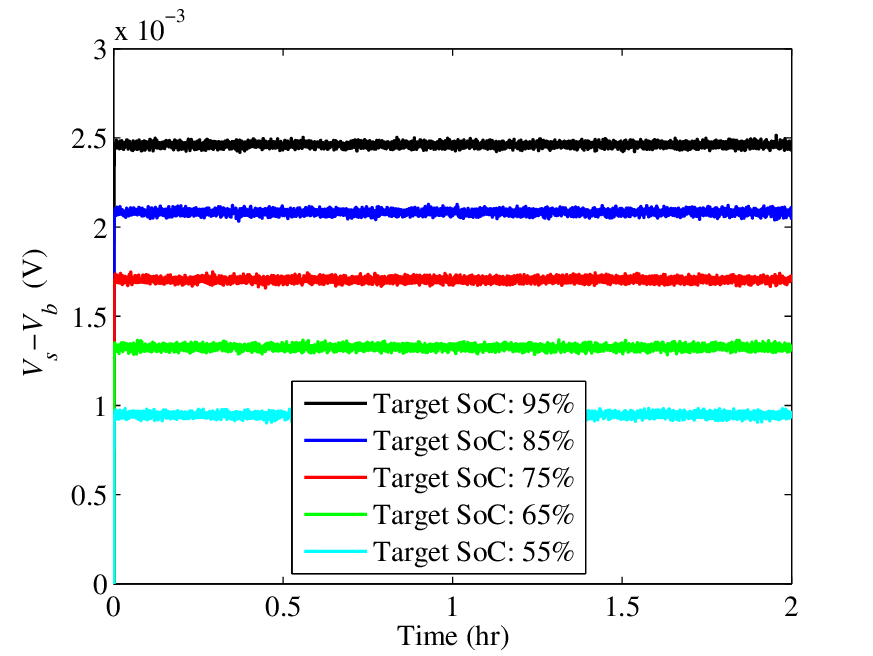}\label{Ex1-CC-Volt-Diff}}
  \caption{Example 1 - Application of \lqcwfts to charge the battery from an initial SoC at $30\%$ to $55\%$, $65\%$, $75\%$, $85\%$ and $95\%$: (a) the SoC trajectories over time; (b) the charging current profiles; (c) the output voltage profiles; (d) the potential differences as health indicator; (e) potential difference due to constant current charging.}
\label{Ex1}
\end{figure*}

{\em Example 1 - Application of \lqcwfts:} Suppose that the user wants to complete the charging in 2 hours. The total number of time instants thus is $N=7200$. Meanwhile, he/she specifies the target \soc value. For the simulation purpose, different target SoC values, $55\%$, $65\%$, $75\%$, $85\%$ and $95\%$, are set here. We apply the \lqcwfts method to carry out the charging tasks. For the control run, $Q_k = 0.1 \cdot (5\times 10^7)^{k/N}$ and $R=0.1$. The exponential increase of $Q_k$ is due to the growing vulnerability of the battery to a larger charging current when the SoC increases. {The practical system will be subject to certain noises, the covariances of which should be included in the Kalman filter implementation. Here, we assume that $W = 10^{-4} I$ and $V = 10^{-4}$.}


\begin{figure*}  \centering
  \vspace{-6mm}
    \subfloat[]{
    \includegraphics[width=0.45\textwidth]{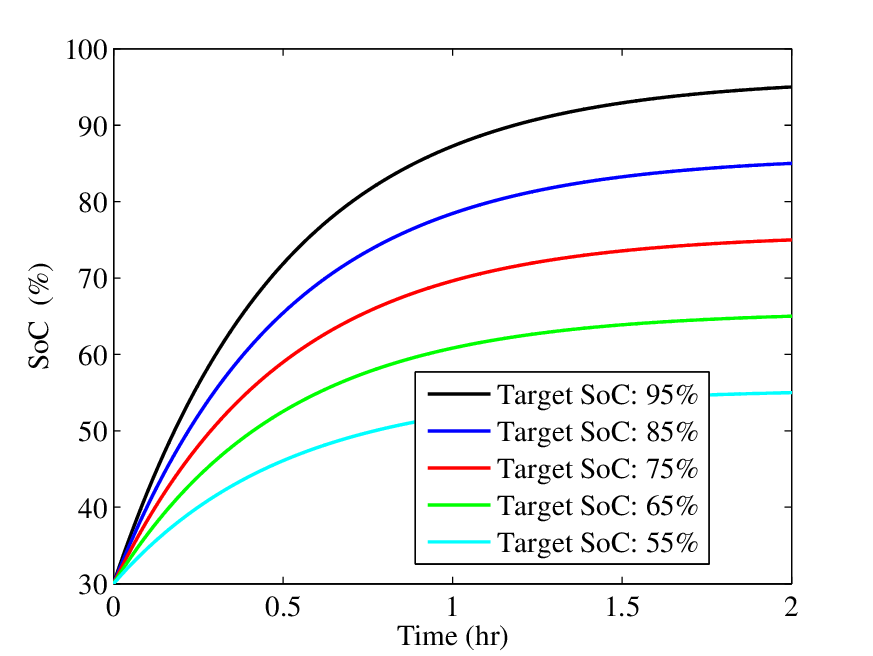}\label{Ex2-SoC}}
  \subfloat[]{
    \includegraphics[width=0.45\textwidth]{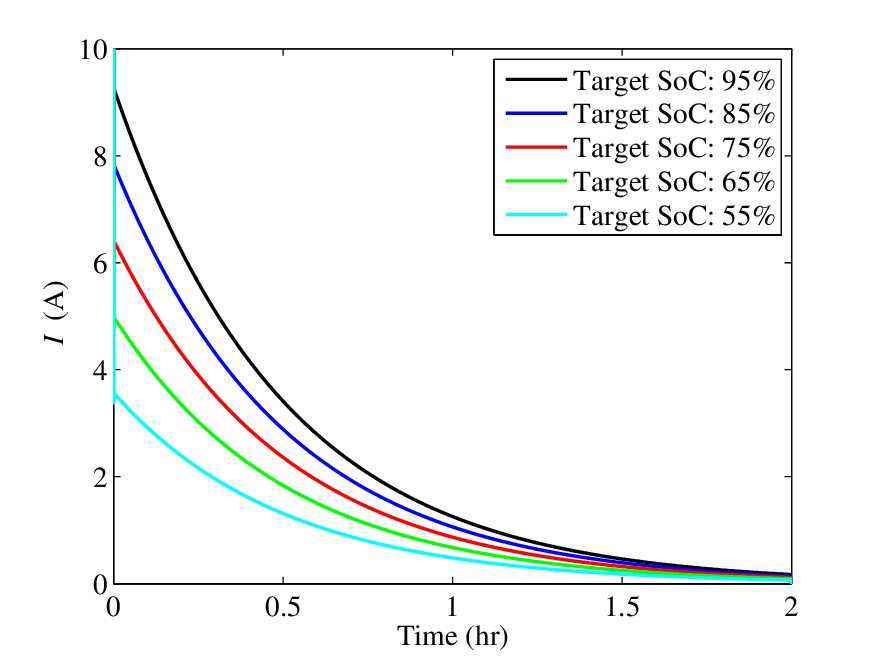}\label{Ex2-Current}}\\
 \subfloat[]{
   \includegraphics[width=0.45\textwidth]{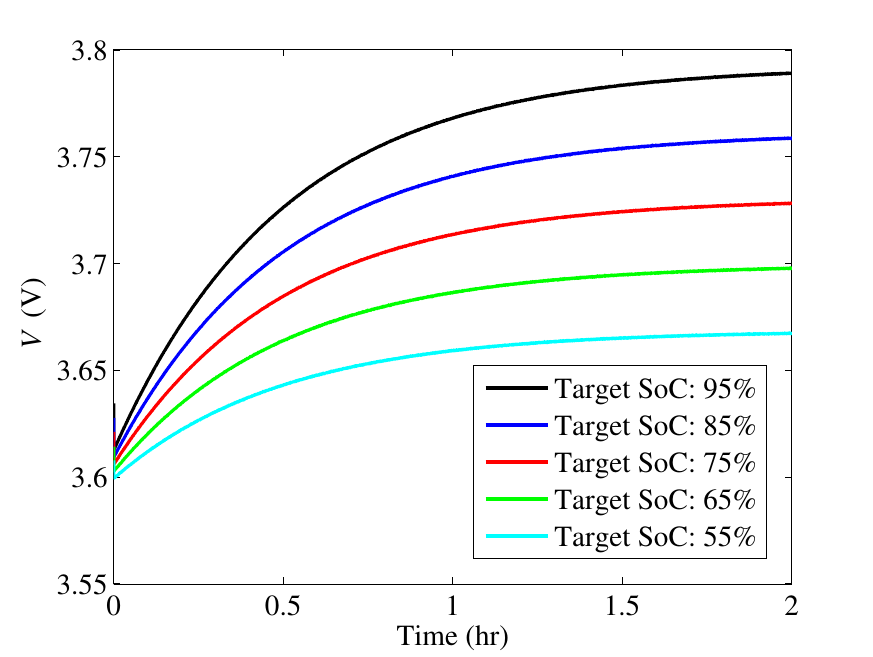}\label{Ex2-Voltage-Output}}
   \subfloat[]{
   \includegraphics[width=0.45\textwidth]{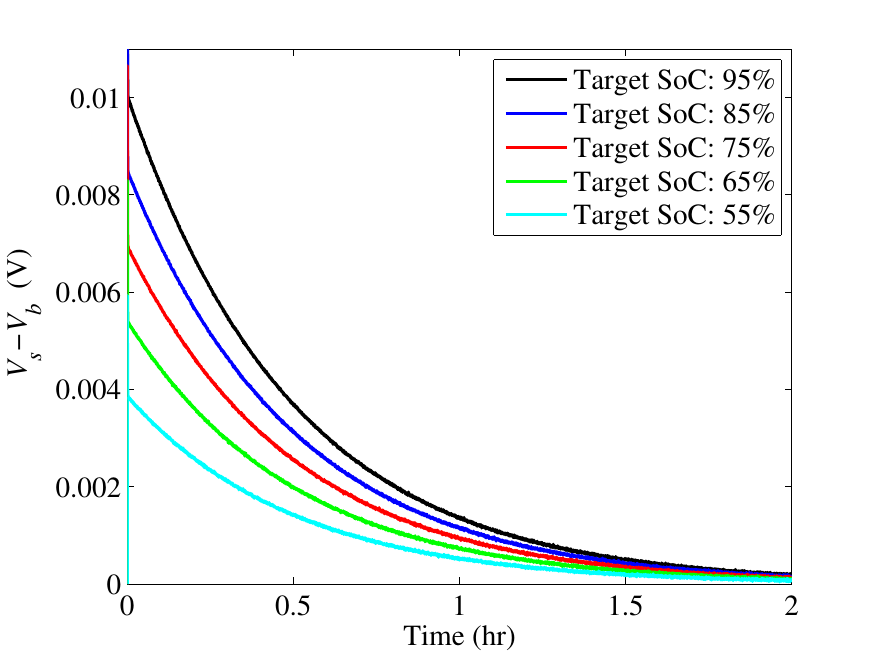}\label{Ex2-Voltage-Diff}}\\
\subfloat[]{
   \includegraphics[width=0.45\textwidth]{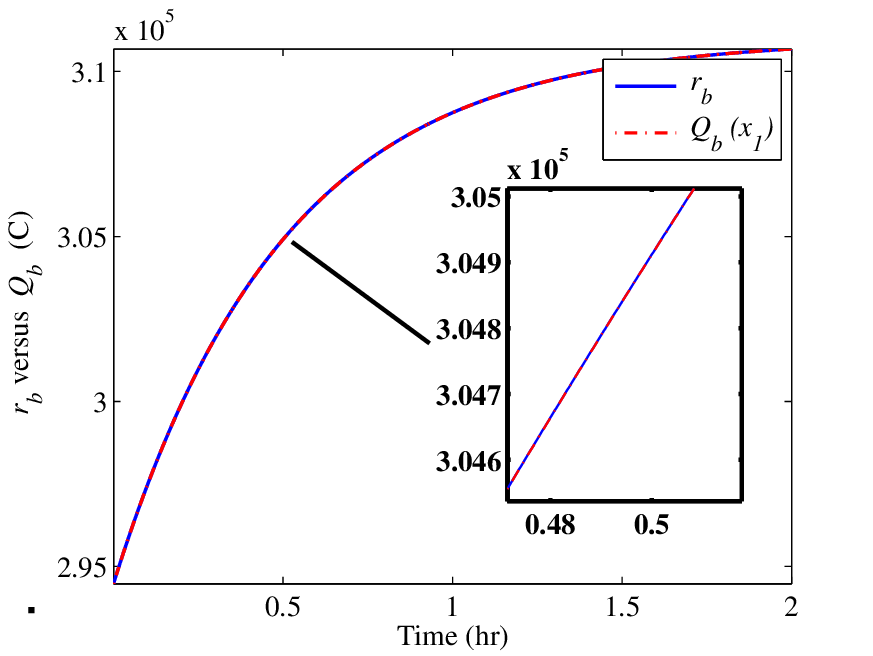}\label{Ex2-Tracking-x1}}
   \subfloat[]{
   \includegraphics[width=0.45\textwidth]{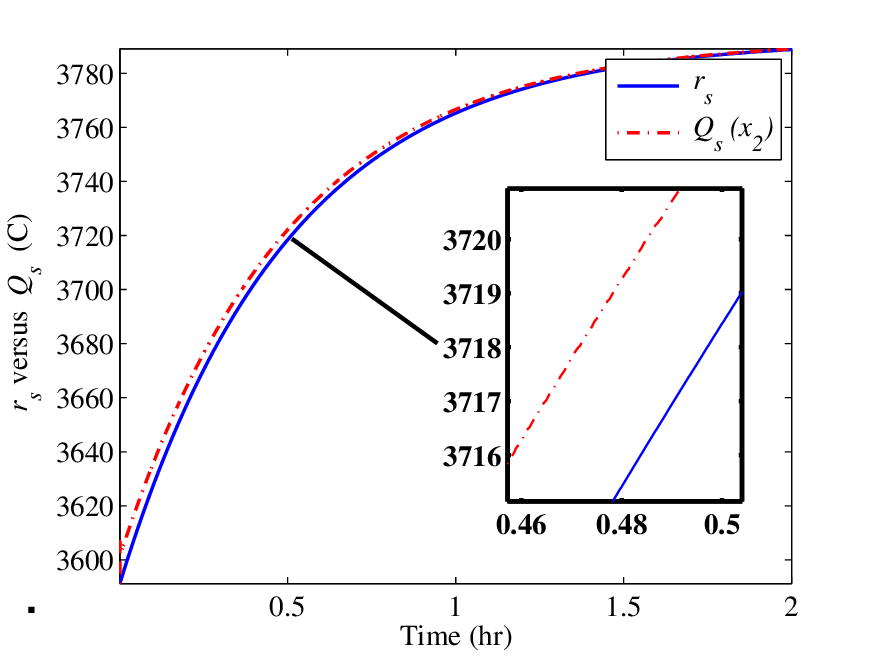}\label{Ex2-Tracking-x2}}
  \caption{Example 2 - Application of \lqt to charge the battery from $30\%$ to $55\%$, $65\%$, $75\%$, $85\%$ and $95\%$: (a) the SoC trajectories over time; (b) the charging current profiles; (c) the output voltage profiles; (d) the potential differences; (e) tracking of $x_1$ (i.e., $Q_b$) for $95\%$ target SoC; (f) tracking of $x_2$ (i.e., $Q_s$) for $95\%$ target SoC.}
\label{Ex2}
\end{figure*}

The computational results are illustrated in Figure~\ref{Ex1}. It is shown in Figure~\ref{Ex1-SoC} that the different target SoCs are satisfied when the charging ends right after two hours, meeting the user-specified objectives. The SoC increases approximately proportionally with time for the first 1.25 hours. Then the rate slows down gradually to zero as the charging objective is being approached. This results from a much larger weight $Q_k$ in the later stage for health protection. The charging current is kept at almost a constant level initially during each charging implementation, as illustrated in Figure~\ref{Ex1-Current}. For a higher target SoC, the magnitude of the current is larger accordingly. However, the current drops quickly in each case as the SoC grows further. {The profiles of the corresponding output voltage is shown in Figure~\ref{Ex1-Volt-Output}. They in general follow a similar trend with the SoC trajectories, rising steadily at first and then at gradually declining rates.}
The voltage difference between $C_s$ and $C_b$, which quantifies the harm incurred to the battery, is characterized in Figure~\ref{Ex1-Voltage-Diff}. In each case, $V_s-V_b$ remains around a constant value in the first hour, despite high-frequency fluctuations due to noise. This is because a battery can accept a higher current at a low \soc level.  Yet the differences decreases drastically as more charge is sent into the battery, in order to maximize the health of the battery's internal structure. {For comparison, we enforce a constant current of appropriate magnitude to flow through the battery for 2 hours to reach the desired SoC. The consequent potential differences are shown in Figure~\ref{Ex1-CC-Volt-Diff}, which are kept at a fixed level unsurprisingly. This, however, will cause much more detrimental effects to the battery when \soc grows, thus expediting the aging processes.}

{\em Example 2 - Application of \sslqt:} We consider the use of \sslqt for charging in this example, which is an upgraded version of \lqt but more computationally efficient. The problem setting and the tasks are the same as in Example 1 --- charging the battery from an SoC of $30\%$ to $55\%$, $65\%$, $75\%$, $85\%$ and $95\%$ in 2 hours for the same battery. The charging trajectory is generated based on the task. For simplicity and convenience, we assume that the desired trajectories for $x_1$ and $x_2$, denoted as $r_b$ and $r_s$, are generated by
\begin{align*}
r_{j,k} = \dfrac{1-e^{-kt_s/\tau_j}}{1-e^{-Nt_s/\tau_j}}(r_{j,N}-r_{j,0})+r_{j,0},
\end{align*}
where $j=b \ {\rm or} \ s, \  k = 1,2,\cdots, N-1$ and $r_{j,0}$ is the initial charge, $r_{j,N}$ the target charge, and $\tau_j$ the time coefficient for $j=b$ or $s$. Note that $r_{j,0}$ and that $r_{j,N}$ can be calculated from the initial \soc and user-specified target \soc. The resultant trajectories have a steep increase followed by a gentle slope, which are reasonable in view of health protection.
Letting $\tau_b = \tau_s = Nt_s/4$, $V_s$ and $V_b$ are forced to be equal through the charging process. Thus at the trajectory design stage, we put the minimization of the detrimental effects well into consideration.

With the reference trajectories generated, the \sslqt strategy is applied to charging. The actual SoC increase over time is demonstrated in Figure~\ref{Ex2-SoC}. All the targets are reached. In each case, the \soc grows at a fast rate when the SoC is at a low level but at a slower rate when the SoC becomes higher. Figure~\ref{Ex2-Current} shows the current produced by \sslqt. The current usually begins with a large magnitude but decreases quickly. {Figure~\ref{Ex2-Voltage-Output} demonstrates the output voltage profiles of the battery, which see a progressively decelerating growth.} The potential difference, given in Figure~\ref{Ex2-Voltage-Diff}, has a similar trend to the current profiles. It is relatively high when the charging starts, and then reduces fast. The state tracking for the task of $95\%$ SoC is shown in~Figures~\ref{Ex2-Tracking-x1} and~\ref{Ex2-Tracking-x2}. It is observed that tracking of $r_b$ by $x_1$ exhibits high accuracy. Tracking of $r_s$ by $x_2$, however, is increasingly accurate, despite a minor deviation in the first hour. Overall, the closer the target SoC is approached, the smaller the tracking error becomes.

In the above examples, different charging current profiles are generated for the same charging task. While the contributory factors include the selection of $Q$ and the reference charging trajectory generation, such a difference poses another important question: how to assess and compare the charging strategies? There is no clear-cut answer yet as it involves a mix of battery electrochemisitry, charging performance, computational complexity, economic cost, and even user satisfaction. Though beyond scope of this paper, evaluation of charging strategies through theoretical analysis and experimental validation will be part of our future quest.

\section{Conclusions}\label{Conclusions}

Effective battery charging management is vital for the development of EVs. Recently, fast charging control has attracted some research effort. However, the problem of health-aware and user-involved charging has not been explored in the literature. In this paper, we propose a set of first-of-its-kind charging strategies, which aim to meet user-defined charging objectives with awareness of the hazards to health. They are developed in the framework of linear quadratic control. One of them is built on control with fixed terminal state, and the other two on tracking a reference charging trajectory. In addition to the merits of health consciousness and user involvement, they are more computationally competitive than most existing charging techniques requiring online real-time optimization solvers. The usefulness of the proposed strategies is evaluated via a simulation study. This work will provide further incentives for research on EV charging management and is also applicable to other battery-powered applications such as consumer electronics devices and renewable energy systems. {Our future research will include battery-type-specific voltage difference limit identification, optimal charging trajectory generation, and a comprehensive assessment of  the charging strategies.}


{

\renewcommand{\theequation}{A.\arabic{equation}}
\setcounter{equation}{0}  
\section*{Appendix}  

\section*{On Approximate Equivalence between~\eqref{RC-model-Eq} and SPM}

This appendix is to present a proof of approximate mathematical equivalence between the RC model in~\eqref{RC-model-Eq} and the SPM. This will demonstrate that the difference in voltages across $C_b$ and $C_s$ in Figure~\ref{RCM-1} approximates the Li-ion concentration gradient in the SPM.

The SPM simplifies each electrode as a spherical particle with area equivalent to the active area
of this electrode~\cite{Santhanagopalan_06}. Striking a balance between mathematical complexity and fidelity toward
capturing key physical and electrochemical phenomena, it has found significant use in the study battery management, e.g.,~\cite{Moura:DSMC:2013,Domenico:DSMC:2010,Fang:CEP:2014}. At the core of the SPM is the 
conservation of Li ions in the electrode phase. Specifically, the migration of
Li ions inside a solid particle is caused by the gradient-induced
diffusion. It follows from the Fick's laws of diffusion that
        \begin{equation}\label{diffusion}
        \frac{\partial c_{j}(r,t)}{\partial t} = \frac{1}{r^2}
        \frac{\partial}{\partial r}\left( D_{j} r^2 \frac{\partial
        c_{j}(r,t)}{\partial r} \right),
        \end{equation}
where $c$ is the concentration of Li ions in the solid electrode, $D$ the diffusion coefficient, $r$ the radial dimension of the spherical particle representing the electrode, and $j=n,p$ with $n$ for the negative electrode and $p$ for the positive one. The associated initial and boundary conditions are given by
        \begin{align}\label{boundary_condition}
        c_{j}(r,0) = c^0, \ \
        \left.\frac{\partial c_{j}}{\partial r}\right|_{r=0}=0, \ \
        \left.\frac{\partial c_{j}}{\partial r}\right|_{r=R_j}=-\frac{1}{D_{s,j}} J_j.
        \end{align}
Here, $J_j$ is the molar flux at the electrode/electrolyte interface of a
single particle. When $j=n$ and $p$, respectively,
        \begin{align}
        J_n(t) = -\frac{I(t)}{FS_n}, \ \
        J_p(t) = \frac{I(t)}{FS_p},
        \end{align}
where $I$ is the charging ($I>0$) or discharging ($I<0$) current, $S$ the surface area, and $R_j$ the radius of the particle.

\begin{figure}[t]  \centering
    \subfloat[]{
    \includegraphics[width=0.22\textwidth]{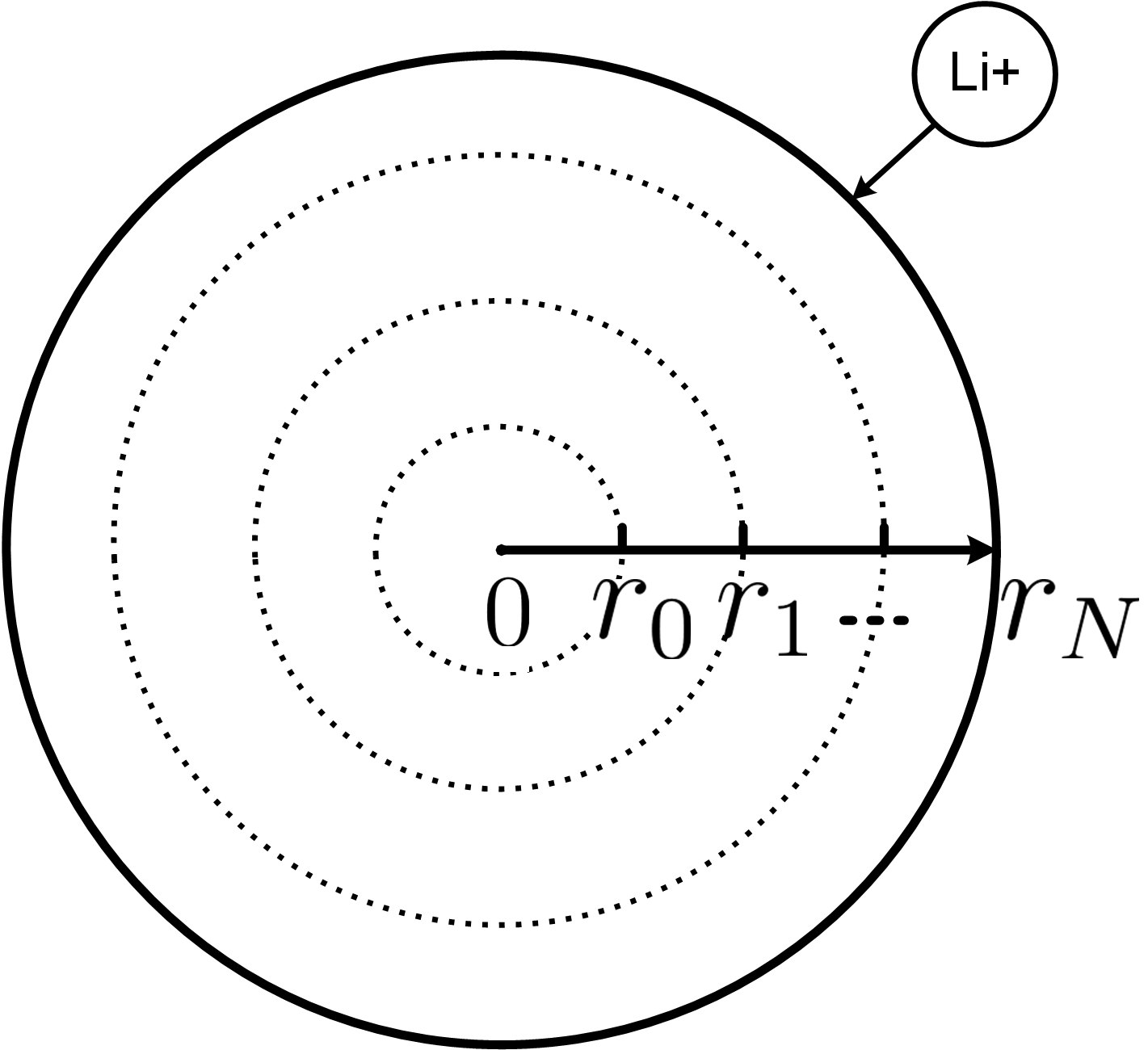}\label{SPM}}\hspace{5mm}
  \subfloat[]{
    \includegraphics[width=0.22\textwidth]{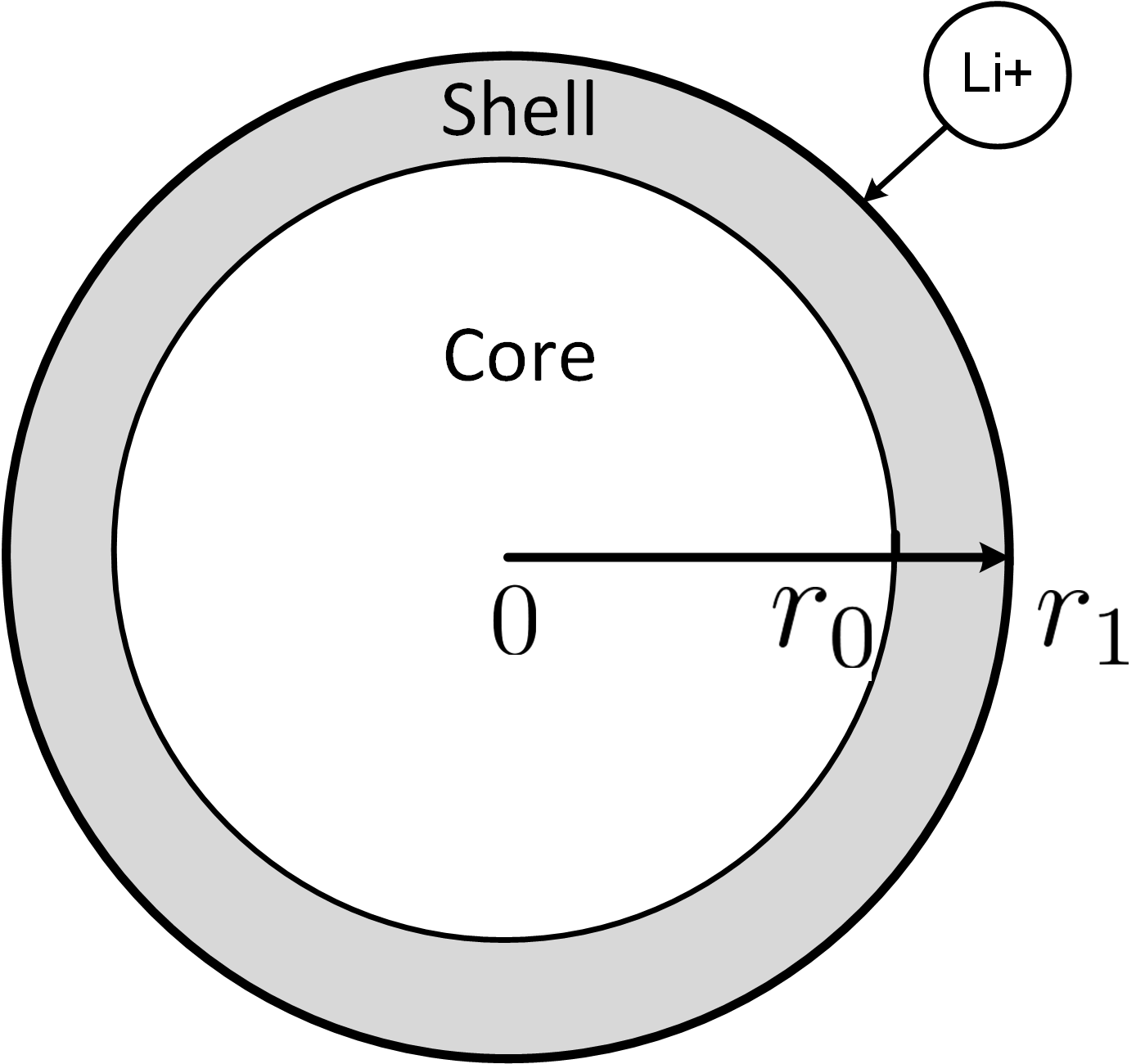}\label{Two-Layer-SPM}} 
  \caption{(a) Subdivision of the spherical particle representing the positive electrode into multiple finite volumes along the radial coordinate; (b) subdivision of the particle into two finite volumes, named the core and shell, respectively.}
\label{Finite-Volume}
\end{figure}


Next, we consider convert the PDE-based diffusion equation into a system of ODE equations using a finite-volume approach. That is, we subdivide the particle along the radial coordinate into a set of continuous finite volumes, as shown in Figure~\ref{SPM}. The finite volume at the center is a ball with a radius $r_0$, and the rest hollow spheres. The $i$-th sphere for $i=0,1,\cdots,N$ has an outer radius of $r_i$ with $r_N = R$. Note that $r_{-1} = 0$.

The total Li-ion amount within the $i$-th finite volume can be quantified as
\begin{align}\nonumber
Q_{j,i} (t) &=\int_{r_{i-1}}^{r_i} c_{j}(r,t)\mathrm{d} V\\ \label{charge_in_FV}
&=  \int_{r_{i-1}}^{r_i} c_{j}(r,t) \cdot 4\pi r^2\mathrm{d} r.
\end{align}
Inserting~\eqref{diffusion} into~\eqref{charge_in_FV}, we have
\begin{align}\label{charge_dynamics} \nonumber
\dot Q_{j,i} (t)
&=  \int_{r_{i-1}}^{r_i} \frac{\partial  c_{j}(r,t)}{\partial t} \cdot 4\pi r^2\mathrm{d} r\\ \nonumber
&= \int_{r_{i-1}}^{r_i} {\rm d} \left(4\pi  D_{j} r^2 \frac{\partial
        c_{j}(r,t)}{\partial r}  \right)\\ \nonumber
&= \left.4\pi  D_{j} r^2 \frac{\partial
        c_{j}(r,t)}{\partial r} \right|_{r_{i-1}}^{r_i}\\
&= -4\pi  D_{j} r_{i-1}^2 \left.\frac{\partial
        c_{j}(r,t)}{\partial r} \right|_{r_{i-1}}+4\pi  D_{j} r_i^2 \left.\frac{\partial
        c_{j}(r,t)}{\partial r} \right|_{r_i}.
\end{align}
To proceed, we assume that the Li ions are uniformly distributed within each finite volume. That is, the Li-ion concentration for the $i$-th sphere is
\begin{align*}
c_j(r,t) = \frac{Q_{j,i}(t)}{ \Delta V_i} \ \mbox{for} \ r_{i-1} < r \leq r_i,
\end{align*}
where $\Delta V_i = 4\pi (r_i^3 -  r_{i-1}^3)/3$. Then, the concentration gradient at $r_i$ can be approximated as
\begin{align}\label{concentration_gradient} \nonumber
 \left.\frac{\partial c_{j}(r,t)}{\partial r} \right|_{r_i} &= \frac{\frac{Q_{j,i+1}(t) }{\Delta V_{i+1}}
- \frac{Q_{j,i}(t) }{\Delta V_i}}{\frac{r_{i+1}-r_{i-1}}{2}} \\ &= \frac{Q_{j,i+1}(t) }{\Delta V_{i+1} \Delta r_{i+1}}
- \frac{Q_{j,i}(t) }{\Delta V_i \Delta r_{i+1}},
\end{align}
where $\Delta r_{i+1} = (r_{i+1} - r_{i-1})/2$.
Then, according to~\eqref{charge_dynamics}-\eqref{concentration_gradient} and the boundary conditions in~\eqref{boundary_condition}, we obtain
\begin{align}\label{charge_dynamics_1}
\dot Q_{j,0} (t) &= - \frac{4 \pi D_j r_0^2 }{\Delta V_0 \Delta r_1} Q_{j,0}(t) + \frac{4 \pi D_j r_0^2}{ \Delta V_1 \Delta r_1 } Q_{j,1}(t),\\ \label{charge_dynamics_2} \nonumber
\dot Q_{j,i} (t) &=  \frac{4 \pi D_j r_{i-1}^2 }{\Delta V_{i-1} \Delta r_i} Q_{j,i-1}(t)\\   & \nonumber \quad - 4 \pi D_j \left(\frac{ (r_{i-1}^2 }{\Delta V_i r_i} + 
\frac{r_i^2}{\Delta V_i \Delta r_{i+1}} \right)Q_{j,i}(t)\\ \nonumber & \quad + \frac{4 \pi D_j r_i^2}{ \Delta V_{i+1} \Delta r_{i+1} } Q_{j,i+1}(t),\\ & \quad \mbox{for} \quad \ 1\leq i < N, \\ \label{charge_dynamics_3} \nonumber
\dot Q_{j,N} (t) & = \frac{4 \pi D_j r_{N-1}^2 }{\Delta V_{N-1} \Delta r_N} Q_{j,N-1}(t) - \frac{4 \pi D_j (r_{N-1}^2 ) }{\Delta V_N \Delta r_N } Q_{j,N}(t)\\& \quad \pm \frac{4\pi r_N^2}{F S_{p(n)}} I(t).
\end{align}

Now let us consider only the positive electrode without loss of generality and suppose that its particle is subdivided into only two finite volumes, the bulk inner domain (core) and the near-surface domain (shell), with $r_0 \gg r_1-r_0$. This approximates the charge diffusion at the interface between the near-surface area and the inside of the particle. By~\eqref{charge_dynamics_1}-\eqref{charge_dynamics_3}, we have
\begin{align}\label{simple_SPM}
\left[
\begin{matrix}
\dot Q_{p,0} (t) \cr \dot Q_{p,1}(t)
\end{matrix}
\right]
=
\left[
\begin{matrix}
- \eta_0 &  \eta_1\cr  \eta_0 & - \eta_1
\end{matrix}
\right] 
\left[
\begin{matrix}
Q_{p,0} (t) \cr  Q_{p,1}(t)
\end{matrix}
\right] + \left[
\begin{matrix}
0 \cr \gamma
\end{matrix}
\right] I(t)
\end{align}
where $\eta_0 = {4 \pi D_p r_0^2 }/{\Delta V_0 \Delta r_1} $, $\eta_1 = {4 \pi D_p r_0^2 }/{\Delta V_1 \Delta r_1} $, and $\gamma = 4\pi r_0^2 /F S_p$. 


It is seen that $\eta_0 \ll \eta_1$ in~\eqref{simple_SPM} due to $\Delta V_0 \gg \Delta V_1$ and that $R_s/(R_b+R_s)$ is close to 0 because $R_s \ll R_b+R_s$ in~\eqref{RC-model-Eq}. With this observation and comparing~\eqref{simple_SPM} with~\eqref{RC-model-Eq}, we can find that they share an approximately equivalent mathematical form. 
Thus, from the perspective of physical abstraction, we can associate the shell of the particle with the surface capacitor $C_s$ and the core with the bulk capacitor $C_b$.  Meanwhile, an analogy can be drawn between the voltage difference $V_s-V_b = Q_s/C_s - Q_b/C_b$ and the gradient of the Li-ion concentration in the two finite volumes, which is expressed as
$Q_{p,1}/\Delta V_1  - Q_{p,0}/\Delta V_0
$. This finding justifies the use of the voltage difference in Sections II and III.
}

\balance

\bibliographystyle{IEEEtran}
\bibliography{LQT_Charging}

\begin{thebibliography}{10}
\providecommand{\url}[1]{#1}
\csname url@samestyle\endcsname
\providecommand{\newblock}{\relax}
\providecommand{\bibinfo}[2]{#2}
\providecommand{\BIBentrySTDinterwordspacing}{\spaceskip=0pt\relax}
\providecommand{\BIBentryALTinterwordstretchfactor}{4}
\providecommand{\BIBentryALTinterwordspacing}{\spaceskip=\fontdimen2\font plus
\BIBentryALTinterwordstretchfactor\fontdimen3\font minus
  \fontdimen4\font\relax}
\providecommand{\BIBforeignlanguage}[2]{{%
\expandafter\ifx\csname l@#1\endcsname\relax
\typeout{** WARNING: IEEEtran.bst: No hyphenation pattern has been}%
\typeout{** loaded for the language `#1'. Using the pattern for}%
\typeout{** the default language instead.}%
\else
\language=\csname l@#1\endcsname
\fi
#2}}
\providecommand{\BIBdecl}{\relax}
\BIBdecl

\bibitem{EDTA:2014}
\BIBentryALTinterwordspacing
{Electric Drive Transportation Association}. (2015, May) Cumulative {U.S.}
  plug-in vehicle sales. [Online]. Available:
  \url{http://www.electricdrive.org/index.php?ht=d/sp/i/20952/pid/20952}
\BIBentrySTDinterwordspacing

\bibitem{Domenico:DSMC:2010}
D.~Domenico, A.~Stefanopoulou, and G.~Fiengo, ``Lithium-ion battery state of
  charge and critical surface charge estimation using an electrochemical
  model-based extended {Kalman} filter,'' \emph{ASME Journal of Dynamic
  Systems, Measurement, and Control}, vol. 132, no.~6, pp.
  061\,302--061\,302--11, 2010.

\bibitem{Fang:CEP:2014}
H.~Fang, Y.~Wang, Z.~Sahinoglu, T.~Wada, and S.~Hara, ``State of charge
  estimation for lithium-ion batteries: An adaptive approach,'' \emph{Control
  Engineering Practice}, vol.~25, pp. 45 -- 54, 2014.

\bibitem{Wang:TCST:2015}
Y.~Wang, H.~Fang, Z.~Sahinoglu, T.~Wada, and S.~Hara, ``Adaptive estimation of
  the state of charge for lithium-ion batteries: Nonlinear geometric observer
  approach,'' \emph{IEEE Transactions on Control Systems Technology}, 2015, in
  press.

\bibitem{Kim:JPS:2015}
T.~Kim, Y.~Wang, H.~Fang, Z.~Sahinoglu, T.~Wada, S.~Hara, and W.~Qiao,
  ``Model-based condition monitoring for lithium-ion batteries,'' \emph{Journal
  of Power Sources}, vol. 295, pp. 16 -- 27, 2015.

\bibitem{Fang:JPS:2014}
H.~Fang, X.~Zhao, Y.~Wang, Z.~Sahinoglu, T.~Wada, S.~Hara, and R.~A.
  de~Callafon, ``Improved adaptive state-of-charge estimation for batteries
  using a multi-model approach,'' \emph{Journal of Power Sources}, vol. 254,
  pp. 258 -- 267, 2014.

\bibitem{Smith:TCST:2010}
K.~A. Smith, C.~D. Rahn, and C.~Y. Wang, ``Model-based electrochemical
  estimation and constraint management for pulse operation of lithium ion
  batteries,'' \emph{IEEE Transactions on Control Systems Technology}, vol.~18,
  no.~3, pp. 654--663, 2010.

\bibitem{Bartlett:TCST:2016}
A.~Bartlett, J.~Marcicki, S.~Onori, G.~Rizzoni, X.~G. Yang, and T.~Miller,
  ``Electrochemical model-based state of charge and capacity estimation for a
  composite electrode lithium-ion battery,'' \emph{IEEE Transactions on Control
  Systems Technology}, vol.~24, no.~2, pp. 384--399, 2016.

\bibitem{Dey:TCST:2015}
S.~Dey, B.~Ayalew, and P.~Pisu, ``Nonlinear robust observers for
  state-of-charge estimation of lithium-ion cells based on a reduced
  electrochemical model,'' \emph{IEEE Transactions on Control Systems
  Technology}, vol.~23, no.~5, pp. 1935--1942, 2015.

\bibitem{Moura:DSMC:2013}
M.~Scott, N.~A. Chaturvedi, and K.~Miroslav, ``Adaptive partial differential
  equation observer for battery state-of-charge/state-of-health estimation via
  an electrochemical model,'' \emph{ASME Journal of Dynamic Systems,
  Measurement, and Control}, vol. 136, no.~1, 2013.

\bibitem{Lin:TCST:2013}
X.~Lin, H.~Perez, J.~Siegel, A.~Stefanopoulou, Y.~Li, R.~Anderson, Y.~Ding, and
  M.~Castanier, ``Online parameterization of lumped thermal dynamics in
  cylindrical lithium ion batteries for core temperature estimation and health
  monitoring,'' \emph{IEEE Transactions on Control Systems Technology},
  vol.~21, no.~5, pp. 1745--1755, 2013.

\bibitem{Suthar:PCCP:2013}
B.~Suthar, V.~Ramadesigan, S.~De, R.~D. Braatz, and V.~R. Subramanian,
  ``Optimal charging profiles for mechanically constrained lithium-ion
  batteries,'' \emph{Physical Chemistry Chemical Physics}, vol.~16, no.~1, pp.
  277--287, 2013.

\bibitem{Spotnitz:JPS:2003}
R.~Spotnitz, ``Simulation of capacity fade in lithium-ion batteries,''
  \emph{Journal of Power Sources}, vol. 113, no.~1, pp. 72 -- 80, 2003.

\bibitem{Bergveld:Springer:2002}
H.~Bergveld, W.~Kruijt, and P.~Notten, \emph{Battery Management Systems: Design
  by Modeling}.\hskip 1em plus 0.5em minus 0.4em\relax Springer, 2002.

\bibitem{Catherino:JPS:2004}
H.~A. Catherino, F.~F. Feres, and F.~Trinidad, ``Sulfation in lead–acid
  batteries,'' \emph{Journal of Power Sources}, vol. 129, no.~1, pp. 113 --
  120, 2004.

\bibitem{Young:Springer:2012}
K.~Young, C.~Wang, L.~Wang, and K.~Strunz, ``Electric vehicle battery
  technologies,'' in \emph{Electric Vehicle Integration into Modern Power
  Networks}, R.~Garcia-Valle and J.~P. Lopes, Eds.\hskip 1em plus 0.5em minus
  0.4em\relax Springer, 2012.

\bibitem{Rahn:Wiley:2013}
C.~D. Rahn and C.-Y. Wang, \emph{Battery Systems Engineering}.\hskip 1em plus
  0.5em minus 0.4em\relax Wiley, 2013.

\bibitem{Wong:JPS:2008}
Y.~Wong, W.~Hurley, and W.~W\, ``Charge regimes for valve-regulated lead-acid
  batteries: Performance overview inclusive of temperature compensation,''
  \emph{Journal of Power Sources}, vol. 183, no.~2, pp. 783--791, 2008.

\bibitem{Lam:JPS:1995}
L.~T. Lam, H.~Ozgun, O.~V. Lim, J.~A. Hamilton, L.~H. Vu, D.~G. Vella, and
  D.~A.~J. Rand, ``Pulsed-current charging of lead/acid batteries --- a
  possible means for overcoming premature capacity loss?'' \emph{Journal of
  Power Sources}, vol.~53, no.~2, pp. 215--228, 1995.

\bibitem{Purushothaman:JES:2006}
B.~K. Purushothaman and U.~Landau, ``Rapid charging of lithium-ion batteries
  using pulsed currents: A theoretical analysis,'' \emph{Journal of The
  Electrochemical Society}, vol. 153, no.~3, pp. A533--A542, 2006.

\bibitem{Aryanfar:JPCL:2014}
A.~Aryanfar, D.~Brooks, B.~V. Merinov, W.~A. Goddard, A.~J. Colussi, and M.~R.
  Hoffmann, ``Dynamics of lithium dendrite growth and inhibition: Pulse
  charging experiments and {Monte Carlo} calculations,'' \emph{The Journal of
  Physical Chemistry Letters}, vol.~5, no.~10, pp. 1721--1726, 2014.

\bibitem{Klein:ACC:2011}
R.~Klein, N.~Chaturvedi, J.~Christensen, J.~Ahmed, R.~Findeisen, and A.~Kojic,
  ``Optimal charging strategies in lithium-ion battery,'' in \emph{Proceedings
  of American Control Conference}, 2011, pp. 382--387.

\bibitem{Yan:Energies:2011}
J.~Yan, G.~Xu, H.~Qian, and Z.~Song, ``Model predictive control-based fast
  charging for vehicular batteries,'' \emph{Energies}, pp. 1178--1196, 2011.

\bibitem{Perez:Mechatronics:2015}
H.~Perez, N.~Shahmohammadhamedani, and S.~Moura, ``Enhanced performance of
  li-ion batteries via modified reference governors and electrochemical
  models,'' \emph{IEEE/ASME Transactions on Mechatronics}, vol.~20, no.~4, pp.
  1511--1520, 2015.

\bibitem{Torchio:ACC:2015}
M.~Torchio, N.~A. Wolff, D.~M. Raimondo, L.~Magni, U.~Krewer, R.~B. Gopaluni,
  J.~A. Paulson, and R.~D. Braatz, ``Real-time model predictive control for the
  optimal charging of a lithium-ion battery,'' in \emph{Proceedings of American
  Control Conference}, 2015, pp. 4536--4541.

\bibitem{Liu:DSCC:2015}
J.~Liu, G.~Li, and H.~K. Fathy, ``Efficient lithium-ion battery model
  predictive control using differential flatness-based pseudospectral
  methods,'' in \emph{Proceedings of ASME 2015 Dynamic Systems and Control
  Conference}, 2015, p. {V001T13A005}.

\bibitem{Wai:IET-PE:2012}
R.~Wai and S.~Jhung, ``Design of energy-saving adaptive fast-charging control
  strategy for {Li-Fe-PO4} battery module,'' \emph{IET Power Electronics},
  vol.~5, no.~9, pp. 1684--1693, 2012.

\bibitem{Wang:AUTO:2012}
T.~Wang and C.~G. Cassandras, ``Optimal control of batteries with fully and
  partially available rechargeability,'' \emph{Automatica}, vol.~48, no.~8, pp.
  1658--1666, 2012.

\bibitem{Lewis:Wiley:2012}
F.~L. Lewis, D.~L. Vrabie, and V.~L. Syrmos, \emph{Optimal Control},
  3rd~ed.\hskip 1em plus 0.5em minus 0.4em\relax Wiley, 2012.

\bibitem{Duncan:IEEE-TAC:1999}
T.~Duncan, L.~Guo, and B.~Pasik-Duncan, ``Adaptive continuous-time linear
  quadratic {Gaussian} control,'' \emph{IEEE Transactions on Automatic
  Control}, vol.~44, no.~9, pp. 1653--1662, 1999.

\bibitem{Duncan:IEEE-TAC:2013}
T.~Duncan, ``Linear-exponential-quadratic gaussian control,'' \emph{IEEE
  Transactions on Automatic Control}, vol.~58, no.~11, pp. 2910--2911, 2013.

\bibitem{Lee:AUTO:2000}
J.~H. Lee, K.~S. Lee, and W.~C. Kim, ``Model-based iterative learning control
  with a quadratic criterion for time-varying linear systems,''
  \emph{Automatica}, vol.~36, no.~5, pp. 641 -- 657, 2000.

\bibitem{Johnson:EVS:2000}
V.~H. Johnson, A.~A. Pesaran, and T.~Sack, ``Temperature-dependent battery
  models for high-power lithium-ion batteries,'' in \emph{Proceedings of 17th
  Electric Vehicle Symposium}, 2000.

\bibitem{Johnson:JPS:2002}
V.~H. Johnson, ``Battery performance models in {ADVISOR},'' \emph{Journal of
  Power Sources}, vol. 110, no.~2, pp. 321--329, 2002.

\bibitem{Sitterly:IEEE-TSE:2011}
M.~Sitterly, L.~Y. Wang, G.~Yin, and C.~Wang, ``Enhanced identification of
  battery models for real-time battery management,'' \emph{IEEE Transactions on
  Sustainable Energy}, vol.~2, no.~3, pp. 300--308, 2011.

\bibitem{Seaman:JPS:2014}
A.~Seaman, T.-S. Dao, and J.~McPhee, ``A survey of mathematics-based
  equivalent-circuit and electrochemical battery models for hybrid and electric
  vehicle simulation,'' \emph{Journal of Power Sources}, vol. 256, pp.
  410--423, 2014.

\bibitem{Yuan:Springer:2010}
X.-Z. Yuan, C.~Song, H.~Wang, and J.~Zhang, ``{EIS} equivalent circuits,'' in
  \emph{Electrochemical Impedance Spectroscopy in {PEM} Fuel Cells}.\hskip 1em
  plus 0.5em minus 0.4em\relax Springer, 2010, pp. 139--192.

\bibitem{Knutsen:2013}
D.~Knutsen and O.~Will{\'e}n, ``A study of electric vehicle charging patterns
  and range anxiety,'' Uppsala University, Tech. Rep., 2013.

\bibitem{Markel:AABC:2006}
T.~Markel and A.~Simpson, ``Plug-in hybrid electric vehicle energy storage
  system design,'' in \emph{Advanced Automotive Battery Conference}, 2006.

\bibitem{Pinsona:JES:2013}
M.~B. Pinsona and M.~Z. Bazant, ``Theory of {SEI} formation in rechargeable
  batteries: {Capacity} fade, accelerated aging and lifetime prediction,''
  \emph{Journal of the Electrochemical Society}, vol. 160, no.~2, pp.
  A243--A250, 2013.

\bibitem{Woodford:MIT:2013}
W.~H. {Woodford IV}, ``Electrochemical shock: {Mechanical} degradation of
  ion-intercalation materials,'' Ph.D. dissertation, Massachusetts Institute of
  Technology, 2013.

\bibitem{Bandhauera:JES:2011}
T.~M. Bandhauera, S.~Garimellaa, and T.~F. Fullerb, ``A critical review of
  thermal issues in lithium-ion batteries,'' \emph{Journal of the
  Electrochemical Society}, vol. 158, no.~3, pp. R1--R25, 2011.

\bibitem{Bryson:Taylor:1975}
A.~E. {Bryson, Jr.} and {Yu-Chi Ho}, \emph{Applied Optimal Control}.\hskip 1em
  plus 0.5em minus 0.4em\relax Taylor \& Francis Group, 1975.

\bibitem{Chaturvedi:CSM:2010}
N.~Chaturvedi, R.~Klein, J.~Christensen, J.~Ahmed, and A.~Kojic, ``Algorithms
  for advanced battery-management systems,'' \emph{IEEE Control Systems},
  vol.~30, no.~3, pp. 49--68, 2010.

\bibitem{Santhanagopalan_06}
S.~Santhanagopalan, Q.~Guo, P.~Ramadass, and R.~E. White, ``Review of models
  for predicting the cycling performance of lithium ion batteries,''
  \emph{Journal of Power Sources}, vol. 156, no.~2, pp. 620--628, 2006.

\end{thebibliography}

\end{document}